\documentclass[onefignum,onetabnum]{siamart171218}

\usepackage{amsfonts}
\usepackage{amsmath}
\usepackage{lmodern}
\usepackage{mathtools}
\usepackage[T1]{fontenc}
\usepackage{graphicx}
\usepackage{epstopdf}
\usepackage{algorithmic}
\usepackage[T1]{fontenc}
\usepackage{tikz,caption}
\usepackage{subcaption}
\usepackage{pgfplots}
\pgfplotsset{compat=newest}
\pgfplotsset{plot coordinates/math parser=false}
\usepackage{url}
\usepackage{color}
\usetikzlibrary{plotmarks} 
\usetikzlibrary{external}
\usepgfplotslibrary{external} 
\tikzset{external/mode=graphics if exists}

\newlength\figureheight
\newlength\figurewidth 

\ifpdf
  \DeclareGraphicsExtensions{.eps,.pdf,.png,.jpg}
\else
  \DeclareGraphicsExtensions{.eps}
\fi

\usepackage{xcolor}

\newcommand{\R}{\mathbb{R}}

\newcommand{\hk}{^{(d)}}
\newcommand{\D}{\,\mathrm{d}}

\numberwithin{theorem}{section}

\newcommand{\TheTitle}{A low-rank tensor method for PDE-constrained optimization with isogeometric analysis} 
\newcommand{\ShortTitle}{Low-rank method for PDE-constrained optimization with IgA} 
\newcommand{\TheAuthors}{Alexandra B\"unger, Sergey Dolgov, and Martin Stoll}

\headers{\ShortTitle}{\TheAuthors}
\title{{\TheTitle}}

\author{Alexandra B\"unger\thanks{Technische Universit\"at Chemnitz, Department of Mathematics, Chair of Scientific Computing, 09107 Chemnitz, Germany, \email{alexandra.buenger@mathematik.tu-chemnitz.de}}
\and 
Sergey Dolgov\thanks{University of Bath, Claverton Down, BA2 7AY, Bath, United Kingdom. \email{s.dolgov@bath.ac.uk}}
  \and
  Martin Stoll\thanks{Technische Universit\"at Chemnitz, Department of Mathematics, Chair of Scientific Computing, 09107 Chemnitz, Germany, \email{martin.stoll@mathematik.tu-chemnitz.de}}
}
\usepackage{amsopn}

\ifpdf
\hypersetup{
  pdftitle={\TheTitle},
  pdfauthor={\TheAuthors}
}
\fi

\def\tblue{\textcolor{black}}

\definecolor{grey}{rgb}{0.5,0.5,0.5}
\definecolor{darkgreen}{rgb}{0,0.55,0}

\def\tgreen{\textcolor{black}}

\begin{document}

\maketitle

\begin{abstract}
Isogeometric analysis (IGA) has become one of the most popular methods for the discretization of partial differential equations motivated by the use of NURBS for geometric representations in industry and science. A crucial challenge lies in the solution of the discretized equations, which we discuss in this talk with a particular focus on PDE-constrained optimization discretized using IGA. The discretization results in a system of large mass and stiffness matrices, which are typically very costly to assemble. To reduce the computation time and storage requirements, low-rank tensor methods have become a promising tool. We present a framework for the assembly of these matrices in low-rank form as the sum of a small number of Kronecker products. For assembly of the smaller matrices only univariate integration is required. The resulting low rank Kronecker product structure of the mass and stiffness matrices can be used to solve a PDE-constrained optimization problem without assembling the actual system matrices. We present a framework which preserves and exploits the low-rank Kronecker product format for both the matrices and the solution.
We use the block AMEn method to efficiently solve the corresponding KKT system of the optimization problem. We show several numerical experiments with 3D geometries to demonstrate that the low-rank assembly and solution drastically reduces the memory demands and computing times, depending on the approximation ranks of the domain.
\end{abstract}


\section{Motivation}

Isogeometric Analysis (IgA) is a relatively new discretization technique to give an approximate solution to a problem posed by a partial differential equation (PDE) on a given domain $\Omega$. It was introduced by Hughes, Cottrell and Bazilevs in 2005 \cite{CAD}.

In Isogeometric Analysis the physical domain $\Omega$ and the solution space for solving a PDE via the Galerkin method \cite{strang} are parameterized by the same spline functions, typically B-splines or NURBS (Non uniform rational B-splines). These basis functions are globally defined and have large overlapping supports depending on their degrees. This leads to a global representation of the physical domain and the discretization of the PDEs has a high computational complexity increasing exponentially with respect to the problem's dimension \cite{Mantzaflaris_space_time}.

Recently, a lot of effort has been made to find strategies to overcome this drawback and efficiently assemble the arising system matrices. Here, a big focus lies on exploiting the tensor product structure of the basis functions and lowering the overall computational cost of the basis function quadrature, e.g via finding new quadrature rules \cite{Hiemstra, Hughes} or sum factorization \cite{Antolin}.

To reduce the complexity of the integration and ultimately reduce the overall computation time and storage requirements, Mantzaflaris et al. \cite{angelos1,angelos2} developed a low rank tensor method, which exploits the tensor structure of the basis functions and separates the variables of the integrals. The arising system matrices can then be represented in a compact manner as a sum of Kronecker products of smaller matrices which are assembled via univariate integration, lifting the curse of dimensionality from the integration. This is accomplished via an interpolation step and a low rank representation of the resulting coefficient tensor. 

For two-dimensional settings the low rank approximation can be easily realized by a singular value decomposition of the coefficient matrix. However, in higher dimensions we need to decompose a higher-order tensor which is more challenging, both computationally and conceptually, as there exist many different decompositions and definitions of ranks in the dimensions greater than two.

In this paper we combine the low rank method of Mantzaflaris et al. with low rank Tensor Train (TT) calculations \cite{osel-tt-2011, DoOs-dmrg-solve-2011}. Exploiting the tensor product nature of the arising interpolation we can calculate a low rank TT approximation without prior assembly of the full coefficient tensor by means of the alternating minimal energy (AMEn) method \cite{amen}. We further utilize this method to ultimately solve a large scale optimal control problem in a compact low rank block format, exploiting the Kronecker product structure of the system.

We consider an optimal control problem with a parabolic PDE constraint of the form
\begin{align}
 \min_u&& \frac{1}{2}\int_0^T\int_\Omega (y - \hat{y})^2 \D x\D t &+ \frac{\beta}{2}\int_0^T\int_\Omega u^2 \D x&& \label{equation:optimization1} \\
\mbox{s.t.}&& y_t - \Delta y &=u && \mbox{ in } [0,T]\times \Omega, \label{equation:optimization2}\\
&&y &= 0 && \mbox{ on } [0,T]\times \partial \Omega, \label{equation:optimization3}
\end{align}
with a desired state $\hat{y}$ and control $u$ on a given domain $\Omega$ parameterized by B-Splines or NURBS, as described later on. The discretization of (\ref{equation:optimization1}) - (\ref{equation:optimization3}) in this paper will be performed by isogeometric analysis and the workhorse is the representation of two bilinear forms, the mass term $a_m$ and the stiffness term $a_s$ in a discretized low rank format.

We will briefly review the basic ingredients for isogeometric analysis to clarify the terminology and notations used throughout the paper in Section \ref{section:basics}.
In Section \ref{section:LowRank} we review the previously mentioned low rank approach of Mantzaflaris et al. and present a way to exploit the tensor product structure to quickly find a low rank approximation using the TT format.
We then show how an optimal control problem of the format (\ref{equation:optimization1}) - (\ref{equation:optimization3}) is discretized using IgA and state the resulting discrete saddle point problem in Section \ref{section:optimization}. The discretization results in a very large linear system and in Section \ref{section:LRoptimization} we exploit the derived low rank representation to solve this system in a compact format making use of the iterative Block AMEn method \cite{bdos-sb-2016,ds-navier-2017}.

\tgreen{The performance of the low rank scheme is illustrated by various examples in Section \ref{section:examples}. First, we show the performance for approximating both mass and stiffness matrices in the low rank form for domains with different ranks. We then use these approximations to solve computationally challenging PDE-constrained optimization problems in a low rank tensor train format.}

\section{Basics for IgA} \label{section:basics}

In isogeometric analysis, a geometry is represented exactly using a set of B-splines or NURBS functions. The same basis functions are then used to build the solution space to solve a PDE on the geometric domain \cite{CAD}. The term \textit{B-spline} is short for basis spline and denotes a special type of recursively defined spline. Every spline function with a chosen degree, smoothness, and domain partition can be uniquely represented as a linear combination of B-splines with the same degree, smoothness, and domain partition \cite{boor}.

A set of B-splines is uniquely defined by its degree and knot vector. Choosing a degree $p$, we define a vector $\xi$, called the open knot vector as $\xi= \{\hat{x}_1, \hdots \hat{x}_{n+p+1}\}$ with
\begin{equation} 
0= \hat{x}_1 = \hdots = \hat{x}_{p+1} < \hat{x}_{p+2} \leq \hdots \leq \hat{x}_n < \hat{x}_{n+1} = \hdots = \hat{x}_{n+p+1} = 1, 
\label{equation:knotVector}
\end{equation}
where the end knots appear $p+1$ times and for all other knots duplicate knots are allowed up to multiplicity $p$. The parameter $n$ determines the number of resulting B-splines $\beta_{i,p}$ with $i=1,\hdots,n$. 

For each knot vector $\xi$ as in (\ref{equation:knotVector}), the according  B-splines $\beta_{i,p}$ of degree $p$, with $i = 1,\hdots, n$, are uniquely defined by the recursion
\begin{align}
\beta_{i,0}(\hat{x}) &= \begin{cases} 1 & \mbox{if } \hat{x}_i \leq \hat{x} < \hat{x}_{i+1}, \\ 0 & \mbox{otherwise}, \end{cases} \\
\beta_{i,j}(\hat{x}) &= \frac{\hat{x}-\hat{x}_i}{\hat{x}_{i+j} - \hat{x}_{i}} \beta_{i, j-1}(\hat{x}) + \frac{ \hat{x}_{i+j+1} - \hat{x}}{\hat{x}_{i+j+1} - \hat{x}_{i+1}} \beta_{i+1, j-1}(\hat{x}),
\end{align}
where $j = 1,2,\hdots, p$ and $i = 1, \hdots, n$. Each resulting B-spline $\beta_{i,p}$ has the local support $[\hat{x}_i,\hat{x}_{i+p+1}]$.

We use $\mathbb{S}_\xi^p$ to denote the spline space spanned by the B-splines with degree $p$ and knot vector $\xi$. To construct a B-spline curve in a $D$-dimensional space, the B-splines of $\mathbb{S}_\xi^p$ are combined with given values, called the \textit{control points}, $C_1, \hdots, C_n \in \R^D$ .

Given a B-spline space $\mathbb{S}_\xi^p$ and $n$ control points $C_i \in \R^D$, the curve $F:\R \rightarrow \R^D$ defined by 
\begin{equation}
F(\hat{x}) = \sum^n_{i=1} C_i \beta_{i,p}(\hat{x}) \label{equation:B-spline}
\end{equation}
is called a B-spline curve of degree $p$. 

By using a B-spline curve as defined in Equation (\ref{equation:B-spline}), conic geometries can not be represented exactly \cite{NURBS}. 
Conic shapes can only be represented by rational functions. Therefore, a generalization of the B-splines was
developed, the so called NURBS (Non-uniform rational B-splines) \cite{piegl}. The term \emph{non-uniform} refers to the fact that NURBS usually are defined by a knot vector with non-uniformly sized knot spans.

NURBS are used in a wide spectrum of computational application, especially in CAD or CGI environments, where it became the standard tool to model any kind of required shape \cite{NURBS2}. By adding weights to the B-spline functions and rationalizing the curve, a NURBS curve is defined as
\begin{equation}
N(\hat{x}) = \sum^n_{i=1} C_i \frac{\beta_i(\hat{x}) w_i}{\sum^n_{j=1} \beta_j(\hat{x}) w_j } .
\end{equation}

To represent arbitrary $D$-dimensional geometries with B-splines or NURBS as necessary in isogeometric analysis, univariate spline spaces are combined to multivariate spaces via tensor product.

Consider $D$ different univariate spline spaces $\mathbb{S}_{\xi_d}^{p_d}$, each having one-dimensional
variables $\hat{x}\hk \in \mathbb{R}$, with $d = 1,  \hdots ,D$. The knot vector and the spline degree of the according univariate space are denoted by $\xi_d$ and $p_d$.

We obtain a $D$-variate tensor product spline space $\mathbb{S}^D =  \mathbb{S}_{\xi_1}^ {p_1} \otimes \hdots \otimes \mathbb{S}_{\xi_D}^{p_D}$ with variables $\hat{x} = (\hat{x}^{(1)}, \hdots, \hat{x}^{(D)})^T$ as a space of piecewise polynomial functions with degree $p = (p_1, \hdots, p_D)$. Its elements are denoted by
\begin{equation}
\beta_\mathbf{i} (\hat{x}) = \prod_{d=1}^D \beta_{i_d}\hk (\hat{x}\hk).
\end{equation}

Given such a basis $\mathbb{S}^D$, we define a B-spline (or NURBS) geometry mapping $G:\hat{\Omega} \rightarrow \Omega$ from the $D$-dimensional unit cube $\hat{\Omega}:=[0,1]^D$ onto an arbitrary geometric shape $\Omega \subset \R^D$ as
\begin{equation} \label{equation:geometryMapping}
G(\hat{x}) = \sum_{\mathbf{i} \in I} C_\mathbf{i} \beta_{\mathbf{i}}(\hat{x}) = C:B(\hat{x}),
\end{equation}
with control points $C_{\mathbf{i}} \in \R^D$ and multi-index $\mathbf{i} \in I =  \{ (i_1, \hdots, i_D) \, | \, i_d = 1,\hdots, n_d ,\, d=1,\hdots, D\}$. Here $C:B(\hat{x})$ denotes the Frobenius product of a tensor $C \in \R^{D \times n_1 \times \hdots \times n_D}$, holding all the control points, and a tensor $B(\hat{x})\in \R^{n_1\times \hdots \times n_D}$ holding all the basis functions of $\mathbb{S}^D$ evaluated in $\hat{x}$ in a suitable order.

Now that we have a spline representation of the geometry $\Omega$, we can use the same spline functions to parameterize the solution space of a PDE on the geometry. For the discretization of the optimal control problem (\ref{equation:optimization1}) - (\ref{equation:optimization3}) we need to look at the bilinear forms of two stationary problems.

The first important bilinear form, given by
\begin{equation}
a_m(u,v) = \langle u, v \rangle_2 = \int_\Omega uv \D x,
\end{equation}
is called the mass term and results as the weak formulation of the boundary value problem $u(x) = f(x)$ in $\Omega$.

The second bilinear form we will consider is called the stiffness term,
\begin{align}
a_s(u,v) &=  -\int_\Omega  (\Delta u) v \D x, \\
&= \int_\Omega \nabla u \cdot \nabla v \D x, \label{equation:stiffness},
\end{align}
and results from the Poisson equation, $- \Delta u(x) = f(x)$ in $\Omega$.

We want to produce approximations to the solutions $u \in H^1(\Omega)$ with discrete functions $u_h \in V_h \subset H^1(\Omega)$ constructed with B-splines. In isogeometric analysis we use the same splines from the geometry mapping (\ref{equation:geometryMapping}) to parameterize the solution space $V_h$,
\begin{equation}
 V_h = \mbox{span}\{\beta_\mathbf{i} \circ G^{-1} \, \, : \, \, \mathbf{i} \in I \},
\end{equation}
with an index set $I$ such that $\beta_{\mathbf{i}}$ are the elements of $\mathbb{S}^D$.
The functions in $V_h$ are linear combinations of the basis functions with coefficients $u_\mathbf{i}$,
\begin{equation} 
u_h = \sum_{\mathbf{i} \in I} u_{\mathbf{i}} ( \beta_{\mathbf{i}} \circ G^{-1}).
\label{equation:coefficientSet}
\end{equation}
The space $V_h$ is now used for the Galerkin discretization of the mass and stiffness terms, resulting in the discrete mass term,
\begin{equation}
a_{m,h}(u_h, v_h) = \int_\Omega u_h(x) v_h(x) \D x = \int_{\hat{\Omega}} \sum_{\mathbf{i} \in I} u_\mathbf{i} \beta_\mathbf{i}(\hat{x}) \sum_{\mathbf{j} \in I} v_\mathbf{j} \beta_\mathbf{j}(\hat{x}) \omega(\hat{x}) \D \hat{x},
\end{equation}
with $\omega(\hat{x}) = | \det \nabla G(\hat{x})| $ and the discrete stiffness term
\begin{equation}
a_{s,h}(u_h,v_h) = \int_{\Omega} (\nabla u_h(x) ) \cdot \nabla v_h(x) \D x = \int_{\hat{\Omega}}(Q(\hat{x}) \sum_{\mathbf{i} \in I} u_\mathbf{i} \nabla \beta_\mathbf{i}(\hat{x})) \cdot\sum_{\mathbf{j} \in I} v_\mathbf{j} \nabla \beta_\mathbf{j}(\hat{x}) \D \hat{x},
\end{equation}
with $Q(\hat{x}) = \big (\nabla G(\hat{x})^{T} \nabla G(\hat{x})\big )^{-1} | \det \nabla G(\hat{x})|$ \cite{angelos1}.

This has to hold for all $v_h$ from the test space $V_h$, hence for all combinations of $v_\mathbf{i}\beta_\mathbf{i}$. However, as the basis functions $\beta_\mathbf{i}$ are linearly independent, it is sufficient if the equation holds for each $\beta_\mathbf{j}$ separately. Thus, we can rewrite the discrete bilinear forms as matrix-vector products $Au$ with the vectorization $u$ of the coefficient set $u_{\mathbf{i}}$, $\mathbf{i}\in \mathcal{I}$,
where $A$ is realized as a mass matrix $M$ with elements
\begin{equation} \label{equation:massMatrix}
M_{\mathbf{i},\mathbf{j}} = \int_{\hat{\Omega}} \beta_\mathbf{i} \beta_\mathbf{j} \omega \D \hat{x}
\end{equation}
or a stiffness matrix $K$ with elements
\begin{equation} \label{equation:stiffnessMatrix}
 K_{\mathbf{i},\mathbf{j}} = \int_{\hat{\Omega}}( Q \nabla \beta_\mathbf{i} )\cdot \nabla \beta_\mathbf{j}  \D \hat{x} = \sum_{k,l=1}^D \int_{\hat{\Omega}} q_{k,l} \frac{\partial}{\partial \hat{x}_l} \beta_\mathbf{i}  \frac{\partial}{\partial \hat{x}_k} \beta_\mathbf{j} \D \hat{x}.
\end{equation}

During the derivation of the mass and stiffness matrices, we did not pay attention to the tensor product structure of $\mathbb{S}^D$. We can either arrange $M$ and $K$ as matrices or as tensors of size $(\mathbf{n},\mathbf{n}) = (n_1,\hdots,n_D,n_1,\hdots,n_D)$. With this tensor notation the mass and stiffness matrices in a multi-dimensional setting are represented in a compact way.

Let $B$be the tensor of order \tgreen{D} and size $\mathbf{n} = (n_1,\hdots,n_D)$ holding every $\beta_{\mathbf{i}}\in \mathbb{S}^D$ of equation (\ref{equation:geometryMapping}). All combinations of elements $\beta_\mathbf{i} \beta_\mathbf{j}$ which make up the integrands of $M$ \tgreen{are included in} the tensor product $B \otimes B$.

With this consideration we write the mass term as a tensor $M$
\begin{equation} \label{equation:massTensor1}
 M = \int_{\hat{\Omega}} \omega B \otimes B \D \hat{x} \, \, \, \in \R^{\mathbf n \times \mathbf n},
\end{equation} 
with elements coming from equation (\ref{equation:massMatrix}). The stiffness term can be treated similarly. With the tensor gradient we can write it as a tensor $K$
\begin{align} \label{equation:stiffnessTensor1}
 K &= \int_{\hat{\Omega}} [Q \cdot( \nabla \otimes B)] \cdot (\nabla \otimes B) \D \hat{x} \, \, \, \in \R^{\mathbf n \times \mathbf n}.\\
 &= \sum_{k,l=1}^D \int_{\hat{\Omega}} q_{k,l} \frac{\partial}{\partial \hat{x}_l} B \otimes \frac{\partial}{\partial \hat{x}_k} B \D \hat{x}
\end{align}
whose elements are of the form (\ref{equation:stiffnessMatrix}).
The associated mass and stiffness matrices are obtained by reordering the indices since the elements of the mass and stiffness tensors match the elements of the matrices. So far this tensor structure was not exploited but we will need it to reduce the complexity of the assembly procedure.

\section{Low-rank IGA} \label{section:LowRank}
Looking at the mass and stiffness matrices (\ref{equation:massMatrix}) and (\ref{equation:stiffnessMatrix}), we see that their entries are the product of univariate B-splines and a $D$-variate weight function, $\omega(\hat{x})$ or $Q(\hat{x})$. The scalar $\omega(\hat{x}) = |\mbox{det } \nabla G(\hat{x})|$ and the matrix $Q(\hat{x}) = (\nabla G(\hat{x})^{-1}) (\nabla G(\hat{x}))^{-T} \omega(t) \in \R^{D\times D}$ are determined by the geometry mapping. As Mantzaflaris et al. suggest in \cite{angelos1},we can approximate these weight functions by some combination of univariate functions via interpolation,
\begin{equation}
 \omega(\hat{x}) \approx \omega_1(\hat{x}^{(1)}) \cdots \omega_D(\hat{x}^{(D)}).
\end{equation}
For the low rank approximation of the mass and stiffness matrix, we then approximate the arising multidimensional integrals as products of univariate integrals. The integrands are separable after interpolating the weight functions. To further reduce the computation time and storage requirements of the mass and stiffness matrix calculation, the resulting interpolating function is approximated with low rank methods giving low rank approximations of the system matrices \cite{angelos2,angelos1}.

To do so, we interpolate the weight functions by a combination of univariate B-splines of higher order, denoted by the spline space $\tilde{\mathbb{S}}^D$ with suitable knot vectors $\tilde{\xi}_d$ and degrees $\tilde{p}_d$ with $d=1,\hdots,D$. The weight function $\omega(\hat{x})$ of the mass matrix is interpolated as
\begin{equation}
 \omega(\hat{x}) \approx \sum_{\mathbf{j}\in \mathcal{J}} W_{\mathbf{j}}\tilde{\beta}_{\mathbf{j}}(\hat{x})=W:\tilde{B}(\hat{x}),
\end{equation}
where $\tilde{\beta}_\mathbf{j}(\hat{x})$ are the elements of the spline space $\tilde{\mathbb{S}}^D$ and $\tilde{B}(\hat{x})$ is the tensor holding all $\beta_\mathbf{j}$ ordered according to the index set $\mathcal{J}$.
The weight tensor $W$ has the same dimension as the spline space $\tilde{\mathbb{S}}^D$, being $(\tilde{n}_1,\hdots,\tilde{n}_D)$, and we get its entries by interpolating the weight function in a sufficient number of points, namely $\tilde{n} = \tilde{n}_1\cdots\tilde{n}_D$.

As the derivatives of a B-spline or NURBS are again B-splines or NURBS, the weight function $\omega(\hat{x}) = |\det \nabla G(\hat{x})|$ is again a B-Spline (NURBS) of degree $Dp+1$ \cite{angelos1}. Thus we can get an exact interpolation if we choose basis functions of degree $Dp+1$.

Now we can construct canonical low rank representations of the weight tensor,
\begin{equation} \label{equation:canonicalLowRank}
 W \approx \sum_{r=1}^R \bigotimes_{d=1}^D w_r\hk =: W_R,
\end{equation}
with $w_r\hk \in \R^{n_d}$. With this we get a low rank representation of the weight function,
\begin{equation} \label{equation:weightLowRank}
 \omega(\hat{x}) \approx W_R : \tilde{B}(\hat{x}) = \sum_{r=1}^R \prod_{d=1}^D w_r\hk \cdot \tilde{\beta}\hk(\hat{x}\hk).
\end{equation}
Here $\tilde{\beta}\hk(\hat{x}\hk) \in \R^{n_d}$ denotes the vector holding all univariate basis functions evaluated in $\hat{x}\hk$, and ``$\cdot$'' is the scalar product.
The entries of the mass matrix can be approximated using this low rank representation and we can calculate each entry as the sum of products of univariate integrals,
\begin{align}
 M_{\mathbf{i},\mathbf{j}} &= \int_{\hat{\Omega}} \prod_{d=1}^D \beta_{i_d}\hk \beta_{j_d}\hk \sum_{r=1}^R\prod_{d=1}^D w_r\hk \cdot \tilde{\beta}\hk \D \hat{x} \\
 &= \sum_{r=1}^R \prod_{d=1}^D  \int_0^1\beta_{i_d}\hk\beta_{j_d}\hk w_r\hk \cdot \tilde{\beta}\hk \D \hat{x}\hk.
\end{align}
With these univariate integrals we define a univariate mass matrix, which depends on some weight function $\omega$, as
\begin{equation} \label{equation:massMatrix1D}
 M\hk(\omega) = \int_0^1 B\hk \otimes B\hk \omega \, \D \hat{x}\hk,
\end{equation}
where $B\hk\in \R^{n_d}$ is the vector holding all $n_d$ univariate B-splines of $\mathbb{S}^{p_d}_{\xi_d}$.
According to the tensor representation in equation (\ref{equation:massTensor1}), we can finally write the mass matrix as a sum of Kronecker products of small univariate mass matrices (\ref{equation:massMatrix1D}) with $\omega = w_r\hk \cdot \tilde{\beta}\hk$,
\begin{equation} \label{equation:massFinal}
 M = \sum_{r=1}^R \bigotimes_{d=1}^D M\hk(w_r\hk \cdot \tilde{\beta}\hk).
\end{equation}

The same procedure can be applied to the weight function of the stiffness matrix $Q(\hat{x})$. Note that $Q(\hat{x}) \in \R^{D\times D}$, thus we have to apply the interpolation to each entry of $Q$. Similarly to (\ref{equation:weightLowRank}), for each entry of $Q$ we get the canonical low rank representation
\begin{equation}
 q_{k,l}(\hat{x}) \approx V_{k,l,R} : \tilde{B}(\hat{x}) = \sum_{r=1}^R \prod_{d=1}^D v_{k,l,r}\hk \cdot \tilde{\beta}\hk (\hat{x}\hk), \quad \mbox{ for all } k,l=1,\hdots,D,
\end{equation}
with $v_{k,l,r}\hk \in \R^{n_d}$.

Using this low rank method, we approximate the entries of the stiffness matrix as
\begin{align} 
 K_{\mathbf{i},\mathbf{j}} &= \sum_{k,l=1}^D \int_{\hat{\Omega}} \Big ( \prod_{d=1}^D \delta(l,d) \beta_{i_d}\hk \delta(k,d) \beta_{j_d}\hk \Big )\sum_{r=1}^R \prod_{d=1}^D v_{k,l,r}\hk \cdot \tilde{\beta}\hk \D \hat{x},\\
 & = \sum_{k,l=1}^D \sum_{r=1}^R \prod_{d=1}^D \int_0^1 \delta(l,d) \beta_{i_d}\hk \delta(k,d)\beta_{j_d}\hk v_{k,l,r}\cdot \hk\tilde{\beta}\hk \D \hat{x}\hk,
\end{align}
where $\mathbf{j}=(j_1,\ldots,j_D)$, and $\delta(k,d)$ denotes the operator \tblue{acting on $f$ as}
\begin{equation}
 \delta(k,d) f = \begin{cases} \frac{\partial f}{\partial \hat{x}_d} &\mbox{ if } k = d, \\ f & \mbox{ otherwise}. \end{cases}
\end{equation}
To get a representation for the stiffness matrix corresponding to the mass matrix representation in (\ref{equation:massFinal}), we define the $D^2$ univariate stiffness matrices dependent on some weight function \tblue{$q^{(d)}(\hat{x}\hk)$} as
\begin{equation}
 K_{k,l}\hk(\tblue{q^{(d)}}) = \int_0^1 \left (\delta(l,d) B \right ) \otimes \left (\delta(k,d) B \right ) \tblue{q^{(d)}} \D \hat{x}\hk, \quad \mbox{ for } k,l=1,\hdots,D.
\end{equation}
With this and \tblue{$q^{(d)} = v_{k,l,r}\hk\cdot \tilde{\beta}\hk$} the final low rank tensor representation of the stiffness matrix is
\begin{equation} \label{equation:stiffnessFinal}
 K = \sum_{k,l=1}^D \sum_{r=1}^R \bigotimes_{d=1}^D K_{k,l}\hk(v_{k,l,r}\hk \cdot \tilde{\beta}\hk).
\end{equation}

Both (\ref{equation:massFinal}) and (\ref{equation:stiffnessFinal}) rely on an efficient low rank representation of $W$ and $V_{k,l}$ and we need suitable strategies to perform this task. For a two dimensional setting $W$ and $V_{k,l}$ are matrices and a singular value decomposition (SVD) can be applied easily to find a low rank representation \cite{angelos2}. We approximate the matrix $W \in \R^{n_1\times n_2}$ as
\begin{equation}
 W = U \Sigma V^T \approx \sum_{r=1}^R u_r \sigma_r v_r^T = \sum_{r=1}^R (u_r \sqrt{\sigma_r}) \otimes (v_r \sqrt{\sigma_r}).
\end{equation}
with $U\in \R^{n_1\times n_1}$, $V \in \R^{n_2\times n_2}$ and $\Sigma \in \R^{n_1 \times n_2}$ is the rectangular matrix holding the sorted singular values $\sigma_i$, $i=1,\hdots,\min(n_1,n_2)$ on its diagonal. The low rank approximation is derived by truncating the $n-R$ smallest singular values and the corresponding rows of $U$ and $V$.

In higher dimensional settings, the decomposition becomes more challenging and different types of low rank tensor approximations can be applied as well as considering only a partial decomposition, e. g. into a univariate and a bivariate integration in 3D settings \cite{angelos3}.

\tgreen{For the low rank tensor approximation of a $D$ dimensional tensor as in Equation \eqref{equation:canonicalLowRank} there exists a multitude of possible approximations, e.g. the higher order singular value decomposition (HOSVD), or a CPD decomposition. But for our purpose the tensor train (TT) decomposition is suited  best with respect to simplicity and robustness and we will proceed with TT in the rest of the paper. }

A tensor $W$ is said to be in TT format, if it can be written as
\begin{align}
W(i_1,\hdots,i_D) = W_1(i_1)\cdots W_D(i_D),
\label{eq:tt}
\end{align}
where $W_d(i_d)$ is an $R_{d-1}\times R_d$ matrix for each fixed $i_d$, $1 \leq i_d \leq n_d$ and $R_0 = R_D = 1$ \cite{osel-tt-2011}.

By rearranging the matrices $W_d(i_d)$ for $i_d=1,\hdots, n_d$ into $D$ tensors of sizes $R_{d-1}\times n_d \times R_d$ we can rewrite the TT-format into a canonical low rank representation as desired in (\ref{equation:canonicalLowRank}),
\begin{equation} \label{equation:TTranks}
 W = \sum_{r_1=1}^{R_1}\cdots \sum_{r_D=1}^{R_D} \bigotimes_{d=1}^D \mbox{vec}(W_d(r_{d-1},:,r_d)).
\end{equation}

To interpolate the weight function, we \tgreen{inherently} need to solve a large system of equations $\omega(\hat{X}) = W:\tilde{B}(\hat{X})$, where $\hat{X}$ denotes the set of $n = n_1\cdots n_D$ interpolation points. This equation can be rewritten into
\begin{equation} \label{equation:weightEquationSystem}
 \mbox{vec}(\omega(\hat{X})) = \underbrace{\left(B^{(1)}(\hat{X}^{(1)}) \otimes \hdots \otimes B^{(D)}(\hat{X}^{(D)}) \right )}_{A} \mbox{vec}(W).
\end{equation}
\tblue{The matrix $A$ in \eqref{equation:weightEquationSystem} can be very large for a direct solution.
However, we can first approximate $\omega(\hat{X})$ in a tensor decomposition, and then use the Kronecker structure of $A$ for an efficient computation of a tensor decomposition of $W$.
Indeed, assuming that we have a TT format
$$
\omega(\hat{X}) = \sum_{r_1=1}^{R_1}\cdots \sum_{r_D=1}^{R_D} \bigotimes_{d=1}^D \mbox{vec}(\omega_d(r_{d-1},:,r_d)),
$$
we can write the TT format \eqref{equation:TTranks} for $W$ in the form
\begin{equation} \label{equation:weightEquationSolve}
W = \sum_{r_1=1}^{R_1}\cdots \sum_{r_D=1}^{R_D} \bigotimes_{d=1}^D \left[B^{(d)}(\hat{X}^{(d)})\right]^{-1} \mbox{vec}(\omega_d(r_{d-1},:,r_d)),
\end{equation}
which requires solving $d$ linear systems of sizes $n_1,\ldots,n_D$, respectively.
The TT approximation for $\omega(\hat{X})$ could be precomputed by the TT-SVD \cite{osel-tt-2011} or the TT-Cross \cite{ot-ttcross-2010} methods. We refrain from doing so as the computational time is rather small for the full $\omega(\hat{X})$.
}
This strategy and the efficient representation of $M$ and $K$ allow us to tackle PDE-constrained optimal control problems next.

\section{A PDE-constrained optimization model problem} \label{section:optimization}
We recall the optimal control problem,
\begin{align}
 \min_{y,u}&& \frac{1}{2}\int_0^T\int_\Omega (y - \hat{y})^2 \D x\D t &+ \frac{\beta}{2}\int_0^T\int_\Omega u^2 \D x \D t && \\
\mbox{s.t.}&& y_t - \Delta y &=u && \mbox{ in } [0,T]\times \Omega, \\
&&y &= 0 && \mbox{ on } [0,T]\times \partial \Omega, \label{equation:zeroBoundary} \\
&& y(t=0,\cdot) &= 0 && \mbox{ on } \Omega
\end{align}
to a desired state $\hat{y}$ with control $u$ on a given geometry $\Omega$ and time frame $[0,T]$. 

We want to solve this by discretizing in both time and space resulting in a large saddle point problem \cite{saddlePoint, FEM}. Using an implicit Euler scheme for the time discretization of the PDE and the rectangle rule lead to the time-discrete problem

\begin{align}\label{equation:timeDiscrete}
 \min_{y,u} && \sum_{k=1}^{N_t} \frac{\tau}{2} \Big ( \int_{\Omega} (y_k-\hat{y}_k)^2 \D x  &+ \beta \int_{\Omega} u_k^2 \D x \Big ) && \\
 \mbox{s.t.} && \frac{y_{k+1} - y_{k}}{\tau} - \Delta y_{k+1} &=  u_{k+1} && \mbox{ in } \Omega, \mbox{ for } k=1,\hdots,N_t, \\
 && y_k &= 0 && \mbox{ on } \partial \Omega, \mbox{ for } k=1,\hdots,N_t,
\end{align}
with the number of time steps $N_t$ corresponding to the time step size $\tau = T/{N_t}$ and continuous solution $y_k$ in each time step $k$.

Using the Galerkin-based spatial discretization as described in Section \ref{section:basics} leads to the discrete quadratic problem
\begin{align}
 \min_{y,u}&& \sum_{k=1}^{N_t}\frac{\tau}{2} \big ( (y_k-\hat{y}_k)^TM(y_k-\hat{y}_k) &+ \beta u_k^TMu_k \big) &&\\
\mbox{s.t.}&& \frac{M y_{k}-My_{k-1}}{\tau} + Ky_{k} &= Mu_{k} &&\mbox{ for } k = 1:N_t,
\end{align}
where $M$ and $K$ are the mass and stiffness matrix of the chosen discretization. Here the zero boundary conditions (\ref{equation:zeroBoundary}) are integrated in $M$ and $K$ by omitting the boundary nodes. Omitting the notation from the time-continuous problem (\ref{equation:timeDiscrete}), the states are collected in the vector $y = [y_1,\hdots,y_{N_t}]^T$ with each state $y_k$ being the vector of the corresponding coefficients (\ref{equation:coefficientSet}) and accordingly for the control parameters $u_k$. Note that $y_k$ and $u_k$ are vectors of appropriate dimensionality.

A local minimum of the discretized problem satisfies the \textit{Karush-Kuhn-Tucker} (KKT) conditions \cite{opti2}. The KKT conditions state that if $(y^*,u^*)$ is a local minimum which satisfies a certain constraint qualification, then there exists a multiplier vector $\lambda^*$ such that the \textit{Lagrangian} of the problem, here
\begin{multline} 
 \mathcal{L}(y,u) = \sum_{k=1}^{N_t}\biggl ( \frac{\tau}{2} \big ( (y_k-\hat{y}_k)^TM(y_k-\hat{y}_k) + \beta u_k^TMu_k \big) \\
 + \lambda_k  \big (M y_{k}-My_{k-1} + \tau Ky_{k}-\tau Mu_{k} \big)\biggr ), 
\end{multline}
has a saddle point in the local minimum $(y^*, u^*, \lambda^*)$,
\begin{equation}
 \nabla \mathcal{L}(y^*,u^*,\lambda^*) = 0.
\end{equation}
For our discrete optimization problem this results in the conditions 
\begin{align}
 0 &= \nabla_{y_k}\mathcal{L}(y,u,\lambda) = \tau M (y_k-\hat{y}_k) - \lambda_{k+1} M + \lambda_k(M +\tau K),  \\
 0 &= \nabla_{u_k} \mathcal{L}(y,u,\lambda) = \tau \beta M u_k - \lambda_k \tau M, \\
 0 &= \nabla_{\lambda_k} \mathcal{L}(y,u,\lambda) =  M (y_k-y_{k-1}) + \tau K y_k-\tau Mu_k,
\end{align}
for $k = 1,\hdots, N_t$.

We can rewrite these equations into an equation system
\begin{equation} \label{equation:KKTSystem}
 \begin{bmatrix} \tau \mathcal{M} & 0 & \mathcal{K}^T \\ 0 &\tau \beta \mathcal{M} & -\tau \mathcal{M} \\ \mathcal{K} & -\tau \mathcal{M} & 0 \end{bmatrix} \begin{bmatrix} y \\ u \\ \lambda \end{bmatrix}  = \begin{bmatrix} \tau \mathcal{M}\hat{y} \\ 0 \\ 0 \end{bmatrix},
\end{equation}
where $\mathcal{M} = \mathcal{I}_{N_t} \otimes M$ and $\mathcal{K} = \mathcal{I}_{N_t} \otimes \tau K  + C \otimes M$, where $\mathcal{I}$ is the $N_t\times N_t$ identity matrix and $C$ is
\begin{equation}
 C = \begin{bmatrix} 1 & 0 & 0 & \hdots & 0 \\ -1 & 1 & 0 & \hdots & 0 \\ 0 & -1 & 1 & \hdots & 0 \\ \vdots & & \ddots & \ddots&\\ 0 & \hdots & 0 & -1 & 1 \end{bmatrix}.
\end{equation}
The resulting equation system is a saddle point problem as described in \cite{saddlePoint,stoll1,stoll2}. 

\section{Low-rank solvers for the PDE-constrained optimization problem} \label{section:LRoptimization}
The saddle point problem (\ref{equation:KKTSystem}) typically becomes very large, depending on the number of time steps and refinement in the spatial discretization. By exploiting the tensor product structure for both the solution and the coefficients from Section \ref{section:LowRank}, we can reduce the problem to smaller linear systems on the elements of individual TT blocks.

As we can represent the low rank mass and stiffness matrices as sums of Kronecker products, we can rewrite
\begin{align}
 \label{eq:M_kron}
 \mathcal{M} &= \mathcal{I}_{N_t} \otimes \left ( \sum_{r=1}^R \bigotimes_{d=1}^D M_r^{(d)} \right ),\\
 \label{eq:K_kron}
 \mathcal{K} &= \mathcal{I}_{N_t} \otimes \left ( \sum_{k,l=1}^D \sum_{r=1}^R \bigotimes_{d=1}^D K_{k,l,r}^{(d)} \right ) + C \otimes \left ( \sum_{r=1}^R \bigotimes_{d=1}^D M_r^{(d)} \right ).
\end{align}
With this, each block of (\ref{equation:KKTSystem}) becomes a sum of Kronecker products of small matrices. This structure can be preserved and exploited in appropriate linear solvers, such as Alternating Linear Scheme (ALS) \cite{holtz-ALS-DMRG-2012}, Density Matrix Renormalization Group \cite{schollwock-2011,jeckelmann-dmrgsolve-2002} and Alternating Minimal Energy (AMEn) \cite{amen}.
However, indefinite matrix of saddle point structure in \eqref{equation:KKTSystem} might yield instabilities in the vanilla versions of these algorithms.
We use an extended Block AMEn method (implemented in \textit{amen\_block\_solve.m} in the TT-Toolbox \cite{tt-toolbox}), which preserves the block structure in \eqref{equation:KKTSystem}, and hence the numerical stability.

This algorithm aims to approximate all solution components $y,u,\lambda$ in the same representation, called Block TT decomposition \cite{dkos-eigb-2014}.
Let us collect $y,u,\lambda$ into a matrix
\begin{equation}
f = \begin{bmatrix}f_1 & f_2 & f_3\end{bmatrix} = \begin{bmatrix}y & u & \lambda\end{bmatrix},
\label{eq:w_yup}
\end{equation}
where the components are referred to as $f_{\ell}$, $\ell=1,2,3$.
The Block TT decomposition incorporates the $\ell$-index into one of the factors: instead of \eqref{eq:tt}, we write
\begin{equation}
 f_{\ell}(i_1,\ldots,i_D) = F_{1}(i_1) \cdots F_{d-1}(i_{d-1}) \cdot \hat F_{d}(i_d,\ell) \cdot F_{d+1}(i_{d+1}) \cdots F_{D}(i_D)
 \label{eq:btt}
\end{equation}
for some $1\le d \le D$.
We can put the component enumerator $\ell$ into an arbitrary TT factor using the singular value decomposition.
Suppose we want to move $\ell$ from the factor $d$ to $d+1$.
Consider $\hat F_{d}$ as an $R_{d-1} n_d \times 3 R_d$ matrix,
$
\hat F_{d}(r_{d-1},i_d;~\ell,r_d)
$
and compute the truncated SVD
$$
\hat F_{d} \approx U \Sigma V^\top.
$$
Now we call $U$ the $d$-th TT factor instead of $\hat F_{d}$, and multiply $\Sigma V^\top$ with the $(d+1)$-th factor,
\begin{align}
\label{eq:svd1} F_d(r_{d-1},i_d,r_d')={}&U(r_{d-1},i_d;~r_d'), \\
\label{eq:svd2} \hat F_{d+1}(r_d',i_{d+1},\ell,r_{d+1})={}&\sum_{r_d=1}^{R_d}\Sigma V^\top (r_d';~\ell,r_d) F_{d+1}(r_d,i_{d+1},r_{d+1}).
\end{align}
We have obtained the same representation as \eqref{eq:btt} with $\ell$ sitting in the $(d+1)$-th factor.
This process can be continued further or reversed in order to place $\ell$ in an arbitrary factor.

A state of the art technique for computing directly the factors of a TT decomposition is the Alternating Linear Scheme \cite{holtz-ALS-DMRG-2012,DoOs-dmrg-solve-2011}.
We can observe that the TT representation is linear with respect to the elements of each factor.
Indeed, introduce the following $n^D \times R_{d-1} n_d R_d$ \emph{frame} matrix:
\begin{equation}
  \begin{split}
  F_{\neq d}(i_1,\ldots,i_D;~r_{d-1},j_d,r_d)
   & = F_{1}(i_1) \cdots F_{d-1}(i_{d-1},r_{d-1}) \\
   & \cdot \delta_{i_d,j_d} \\
   & \cdot F_{d+1}(r_d,i_{d+1}) \cdots F_{D}(i_D),
   \end{split}
 \label{eq:frame}
\end{equation}
where $\delta_{i_d,j_d}$ is the identity matrix with respect to the indices $i_d,j_d$.
In case of the Block TT decomposition \eqref{eq:btt}, we assume that we choose the same $d$ for both the position of $\ell$ in \eqref{eq:btt} and the position of the identity matrix in \eqref{eq:frame}.
We can then observe that
\begin{equation}
f_{\ell} = F_{\neq d} \cdot \mbox{vec} (\hat F_{d}(\ell)).
\label{eq:ttlin}
\end{equation}
This linearity allows us to project the original problem into a subspace spanned by the columns of $F_{\neq d}$.
Iterating over all $d=1,\ldots,D$, we obtain the ALS algorithm.
This method starts from some initial guess in the low-rank TT representation,
and hence it never encounters the original (prohibitively large) tensors.

The Block AMEn method \cite{bdos-sb-2016,ds-navier-2017}
projects each of the \emph{submatrices} of \eqref{equation:KKTSystem} onto the frame matrix individually.
For each selected $d=1,\ldots,D$, we compute the elements of $\hat F_{d}$ from the following reduced KKT system:
\begin{equation} \label{eq:KKT_red}
 \begin{bmatrix} \tau F_{\neq d}^T \mathcal{M}F_{\neq d} & 0 & F_{\neq d}^T \mathcal{K}^T F_{\neq d} \\ 0 &\tau \beta F_{\neq d}^T \mathcal{M}F_{\neq d} & -\tau F_{\neq d}^T\mathcal{M}F_{\neq d} \\ F_{\neq d}^T\mathcal{K}F_{\neq d} & -\tau F_{\neq d}^T\mathcal{M}F_{\neq d} & 0 \end{bmatrix}
 \begin{bmatrix} \mbox{vec}~\hat F_{d}(1) \\ \mbox{vec}~\hat F_{d}(2) \\ \mbox{vec}~\hat F_{d}(3) \end{bmatrix}  =
 \begin{bmatrix} \tau F_{\neq d}^T \mathcal{M}\hat{y} \\ 0 \\ 0 \end{bmatrix}.
\end{equation}
This system is small (each submatrix is now of size $R_{d-1} n_d R_d$),
and can be solved efficiently by e.g. MINRES.
Moreover, since $F_{\neq d}$ inherits the TT decomposition of $f$,
and the system matrices $\mathcal{M}$ and $\mathcal{K}$
have the tensor product structure \eqref{eq:M_kron}, \eqref{eq:K_kron},
the reduced matrices $F_{\neq d}^T \mathcal{M}F_{\neq d}$ and $F_{\neq d}^T \mathcal{K}F_{\neq d}$
can be assembled efficiently using the multiplication of tensor trains factor by factor \cite{osel-tt-2011,holtz-ALS-DMRG-2012}.
Having solved \eqref{eq:KKT_red}, we plug the new factor $\hat F_{d}$ back into the block TT decomposition \eqref{eq:btt},
which provides an updated approximation to $y,u$ and $\lambda$ through \eqref{eq:w_yup}.
In order to prepare the next ALS step we move the enumerator $\ell$ to the next factor using SVD (e.g. \eqref{eq:svd1}--\eqref{eq:svd2} in the forward sweep $d\rightarrow d+1$), and recompute the corresponding frame matrix using the new $F_{d}$ factor.
The singular value decomposition makes the appropriate matricization of $F_{d}$ orthogonal, such that the whole frame matrix $F_{\neq d}$ is orthogonal in each step.
This ensures invertibility of the projected matrix in \eqref{eq:KKT_red}.

\section{Numerical experiments} \label{section:examples}

The performance of the low rank tensor train method highly depends on the geometry as the interpolation becomes more challenging and the ranks grow with increasing complexity of the geometry. We conduct some numerical experiments of different complexity to show the advantages of our method compared to the full assembly of the stiffness matrix before combining the assembly with an optimal control problem to show the the performance for large scale saddle point problems.

For our numerical experiments we used \textsc{Matlab} R2018b with the TT-Toolbox \cite{tt-toolbox} on a desktop computer with an \textsc{Intel} Core i7-4770 Quad-Core processor running at $4\times3400$ MHz with 32 GB of RAM.

\subsection{Assembly with TT method}

First, we  will assemble stiffness matrices for different geometries in the low rank format and compare the assembly with a full standard assembly performed by the isogeometric analysis toolbox GeoPDEs 3.0 \cite{geoPDEs} in \textsc{Matlab}. The assembly is compared for different levels of refinement and for each refinement we insert 4 additional knots per knot section in each spatial dimension. We use the same Gauss-Legendre quadrature rule with five quadrature nodes in each spatial direction for both assemblies. For the mass matrix we can get an exact interpolation using a spline space of degree $2p+1$ with $p$ being the degree of the original splines \cite{angelos2}. We use the same degree for the stiffness matrix assembly in our experiments. However, note that we can always increase the degree of the interpolating splines to get a higher accuracy if desired.

\tgreen{Solving \eqref{equation:weightEquationSolve} results in a low rank solution for the given domain and the required truncation tolerance. Hence, the scheme exhibits low-ranks for simpler domains, such as the domain considered next.} The first domain is a three dimensional quarter annulus as shown in Figure \ref{figure:geometry1}. The weight function $Q$ of this geometry is very simple and can be approximated by only one combination of basis functions, thus giving us a rank of 1 for each entry of $Q$.
\tblue{The TT-SVD detects this low rank without prior knowledge of the low rank nature of the geometry, so we can assemble} the stiffness matrix from only one combination of univariate stiffness matrices (\ref{equation:stiffnessFinal}).
We compare the assembly with the full assembly as performed by the geoPDEs toolbox \cite{geoPDEs} and additionally compare the method with an assembly using the CPD method implemented in \cite{tensorlab}. Note that we have to specify the desired rank of the system beforehand to use CPD. Here we chose rank 1, due to the prior knowledge of the low rank nature of the geometry. We use the same quadrature rule for all three assemblies.

 \input{exp_TTassembly.tex}

Figure \ref{figure:assembly1_Time} shows that the TT method is a lot faster than the full assembly and the CPD method, especially for a high number of h-refinements which corresponds to a high number of basis elements. Here the advantage of the TT method over the CPD method lies in the solution of equation (\ref{equation:weightEquationSystem}), which can be efficiently done avoiding full assembly \tgreen{by exploiting the Kronecker product structure} yielding a result already in the low rank Tensor Train format. For the CPD method on the other hand we solve (\ref{equation:weightEquationSystem}) to get a full tensor and \tblue{compute} a CPD of this fully assembled tensor \tblue{using the ALS method. A cross approximation method for CPD could reduce the timings in Figure \ref{figure:assembly1_Time}, but it is significantly less developed and understood than the TT cross scheme, and we are not testing it here.}

Figure \ref{figure:assembly1} further shows that both the CPD and the Tensor Train method approximate the fully assembled stiffness matrix very well. The graph in Figure \ref{figure:assembly1_Diff} depicts the mean difference of the matrices in the Frobenius norm,
\begin{equation}
 \mbox{diff} = \frac{\| S - \tilde{S} \|_F}{\|S\|_F}.
\end{equation}
Note that the CPD method and the TT method find the same decompositions here.
Figure \ref{figure:assembly1_Storage} shows the advantages in regard of storage requirements. For the low rank method we only need to store a small number of small sparse matrices instead of storing the large stiffness matrix reducing the required memory drastically and if a high refinement is desired this effect is amplified.
 
\begin{table}
 \centering
 \begin{tabular}{|c|c|c|c|c|c|c|c|c|}
 \hline
 tolerance & TT-ranks & $q_{1,1}$ & $q_{1,2}$ & $q_{1,3}$ & $q_{2,2}$ & $q_{2,3}$ & $q_{3,3}$ & $\omega$\\
 \hline
 $10^{-10}$ & $R_1$ & 2 & 2 & 1 & 2 & 1 & 11 & 2\\
 & $R_2$ & 2 & 1 & 1 & 1 & 1 & 3 & 1\\
 \hline
 $10^{-7}$ & $R_1$ & 3 & 3 & 1 & 2 & 1 & 7 & 2\\
 & $R_2$ & 1 & 1 & 1 & 1 & 1 & 3 & 1\\
 \hline
$10^{-4}$ & $R_1$ & 2 & 2 & 1 & 2 & 1 & 5 & 2\\
 & $R_2$ & 1 & 1 & 1 & 1 & 1 & 2 & 1\\
\hline
 \end{tabular}
\caption{Ranks for the weight function tensor approximation of figure \ref{figure:geometry3}} \label{table:ranks2}
\end{table}
 
\begin{table}
 \centering
 \begin{tabular}{|c|c|c|c|c|c|c|c|c|}
 \hline
 tolerance & TT-ranks & $q_{1,1}$ & $q_{1,2}$ & $q_{1,3}$ & $q_{2,2}$ & $q_{2,3}$ & $q_{3,3}$ & $\omega$\\
 \hline
 $10^{-7}$ & $R_1$ & 16 & 23 & 24 & 21 & 29 & 21 & 17\\
 & $R_2$ & 9 & 13 &  12 & 11 & 14 & 9 & 9\\ 
 \hline
$10^{-4}$ & $R_1$ & 8 & 10 & 12 & 6 & 15 & 8 & 8\\
& $R_2$ & 5 & 7 & 8 & 8 & 9 & 5 & 5\\
\hline
 \end{tabular}
\caption{Ranks for the weight function tensor approximation of figure \ref{figure:geometry2}} \label{table:ranks1}
\end{table}

 \input{exp_TTassemblyFlag.tex}
 
\input{exp_TTassembly2.tex}

Another example for a low rank domain is the deformed cuboid in Figure \ref{figure:geometry3}. This geometry still possesses a low rank structure and the ranks for different desired accuracies are displayed in Table \ref{table:ranks2}. These ranks stay constant throughout various levels of refinements. We see the comparison of the full assembly and the TT method with different accuracies being $10^{-10}$, $10^{-7}$ and $10^{-4}$ in Figure \ref{figure:assembly3}. We reach the desired accuracies quickly after some refinement steps as depicted in Figure \ref{figure:assembly3_Diff}. Again, the TT method is faster than the full assembly and has an advantage with respect to the storage requirements especially for high refinements, as seen in Figures \ref{figure:assembly3_Time} and \ref{figure:assembly3_Storage}. Note that for the next refinement the fully assembled stiffness matrix would not fit into the memory of our desktop PC anymore.

The performance of the method is remarkable not only for low rank structures as the geometries in Figure \ref{figure:geometry1} or Figure \ref{figure:geometry3} but also for more complex geometries like the high rank domain in Figure \ref{figure:geometry2}. This geometry does not possess a low rank structure but we can still apply the TT method with a high rank or truncate with a desired accuracy to get a low rank approximation.

In Figure \ref{figure:assembly2} we see the comparison of the full assembly and the TT method with different desired accuracies being $10^{-7}$ and $10^{-4}$. The corresponding ranks for the entries of $Q$ are displayed in Table \ref{table:ranks1}. The ranks stay stable and only occasionally vary by $\pm 1$ from the values in table \ref{table:ranks1} due to numerical inaccuracies throughout the different refinement steps.

Recall that the given ranks correspond to the tensor ranks in the TT format (\ref{equation:TTranks}), thus the number of smaller stiffness matrices in each dimension corresponds to the product $R = R_1R_2$. Even though this results in a large number of small matrices, the TT low rank method is still much faster than the full assembly for high refinements.

\subsection{Optimal control examples}

 \input{exp_optimization_rank1.tex}


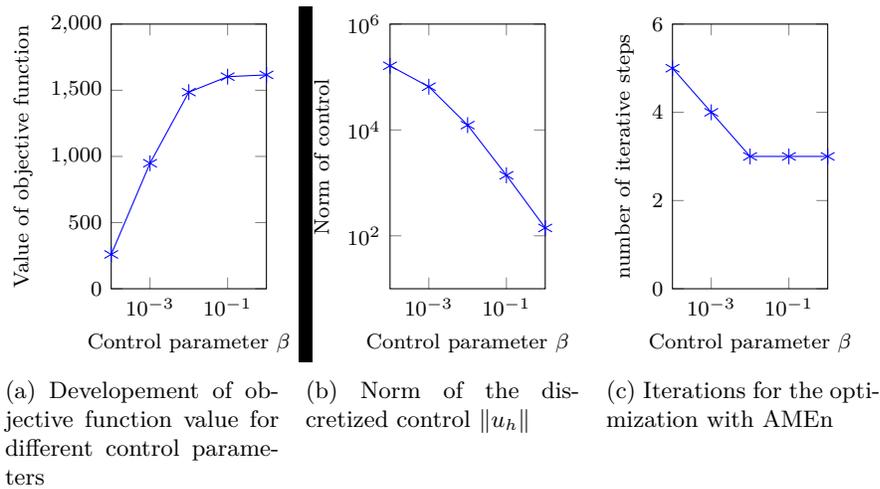
\begin{figure}
 \begin{subfigure}[t]{0.28\textwidth}
  \begin{tikzpicture}[font=\footnotesize]

\begin{axis}[%
width=1\textwidth,
height=1.4\textwidth,
xmode=log,
xmin=0.0001,
xmax=1,
xminorticks=false,
xlabel={$\text{Control parameter }\beta$},
ymin=0,
ymax=2000,
ylabel={Value of objective function},
axis background/.style={fill=white},
]
\addplot [color=blue, mark size=3.0pt, mark=asterisk, mark options={solid, blue}]
  table[row sep=crcr]{%
    0.0001    259.373015595619\\
     0.001    948.462727162855\\
      0.01    1485.69302827373\\
       0.1    1602.74196338784\\
         1    1616.08502958634 \\
};
\end{axis}
\end{tikzpicture}
\caption{Developement of objective function value for different control parameters} \label{figure:rank1opti1_value}
 \end{subfigure}%
 \hspace{1em}%
 \begin{subfigure}[t]{0.28\textwidth}
  \begin{tikzpicture}[font=\footnotesize]
  
\begin{axis}[%
width=1\textwidth,
height=1.4\textwidth,
xmode=log,
xmin=0.0001,
xmax=1,
xminorticks=false,
xlabel={$\text{Control parameter }\beta$},
ymode=log,
ymin=10,
ymax=1000000,
yminorticks=false,
ylabel={Norm of control},
axis background/.style={fill=white},
]
\addplot [color=blue, mark size=3.0pt, mark=asterisk, mark options={solid, blue}]
  table[row sep=crcr]{%
    0.0001    163065.042763983\\
     0.001    65613.4517527527\\
      0.01     12337.540784512\\
       0.1    1389.59170905373\\
         1    140.801849589505\\
};
\end{axis}
\end{tikzpicture}
\caption{Norm of the discretized control $\|u_h\|$} \label{figure:rank1opti_control}
 \end{subfigure}%
 \hspace{1em}%
 \begin{subfigure}[t]{0.28\textwidth}
 \begin{tikzpicture}[font=\footnotesize]
 \begin{axis}[%
width=1\textwidth,
height=1.4\textwidth,
xmode=log,
xmin=0.0001,
xmax=1,
xminorticks=true,
xlabel={$\text{Control parameter }\beta$},
ymin=0,
ymax=6,
ylabel={number of iterative steps},
axis background/.style={fill=white},
]
\addplot [color=blue, mark size=3.0pt, mark=asterisk, mark options={solid, blue}]
  table[row sep=crcr]{%
    0.0001    5 \\
     0.001    4 \\
      0.01    3 \\
       0.1    3 \\
         1    3 \\
};
\end{axis}
\end{tikzpicture}%
\caption{Iterations for the optimization with AMEn} \label{figure:rank1opti_iterations}
  
 \end{subfigure}

 \caption{Performance of the low rank method on the quarter annulus domain from Figure \ref{figure:geometry1} for different control parameters} \label{figure:rank1opti1}

\end{figure}

 \begin{figure}
 \centering
 \begin{subfigure}[t]{0.3\textwidth}
\centering
  \includegraphics[width=0.95\textwidth]{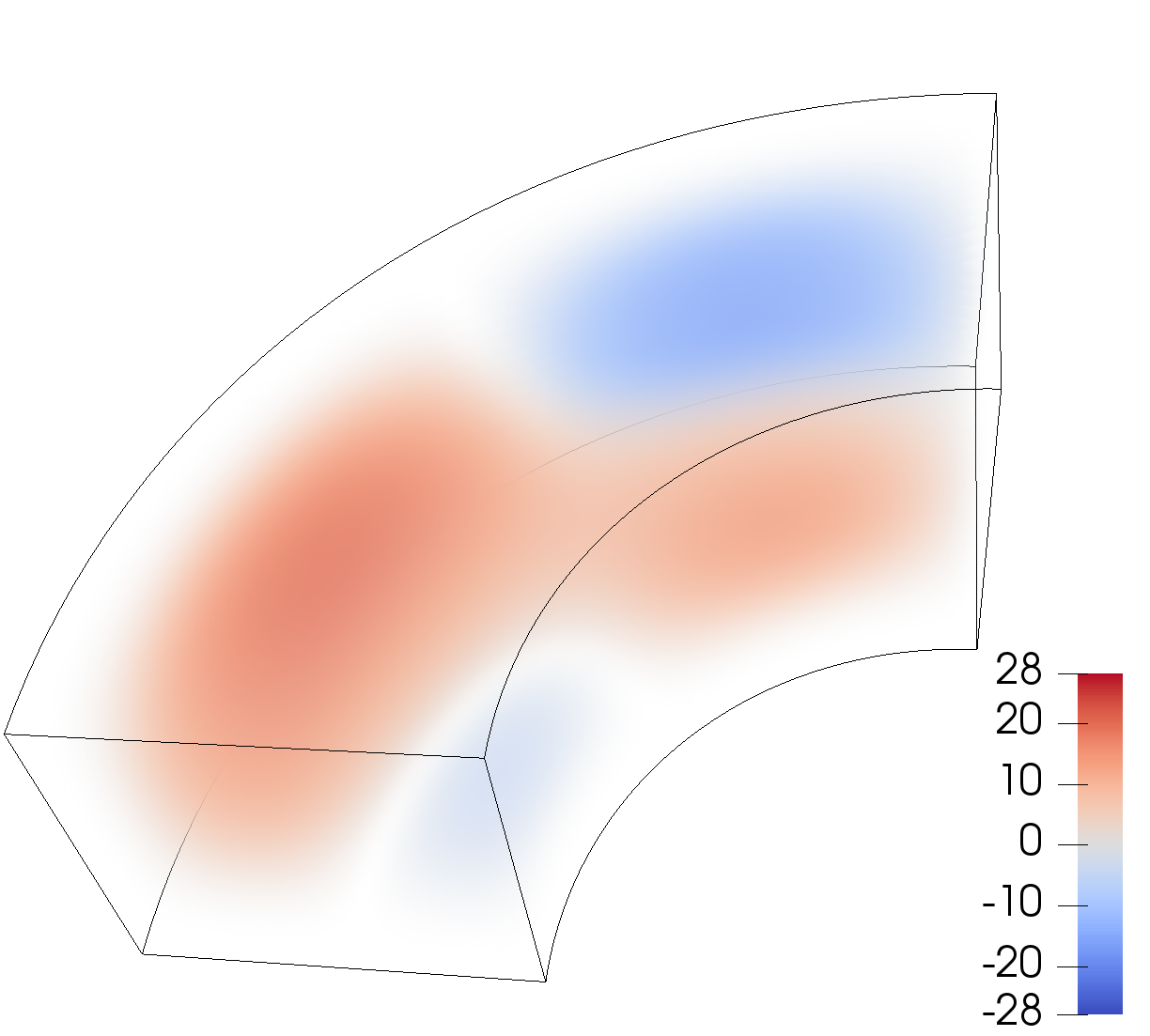}
  \caption{Constant in time desired state $\hat{y}$} \label{figure:rank1_solution_yhat}
  \end{subfigure}%
  \begin{subfigure}[t]{0.3\textwidth}
   \centering
   \includegraphics[width=0.95\textwidth]{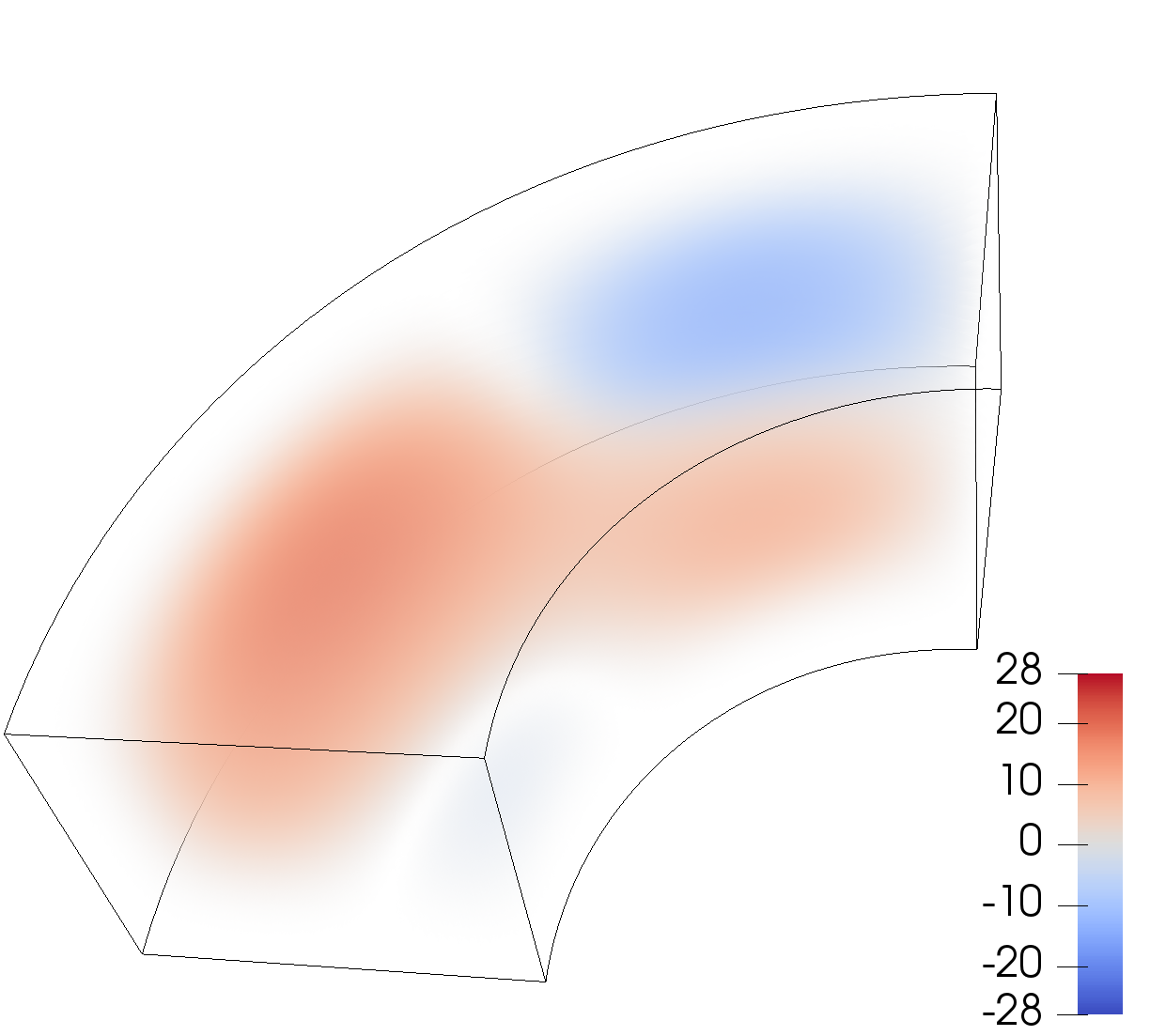}
   \caption{$\beta=10^{-4}$} \label{figure:rank1_solution_1e-4}
  \end{subfigure}%
  \begin{subfigure}[t]{0.3\textwidth}
   \centering
   \includegraphics[width=0.95\textwidth]{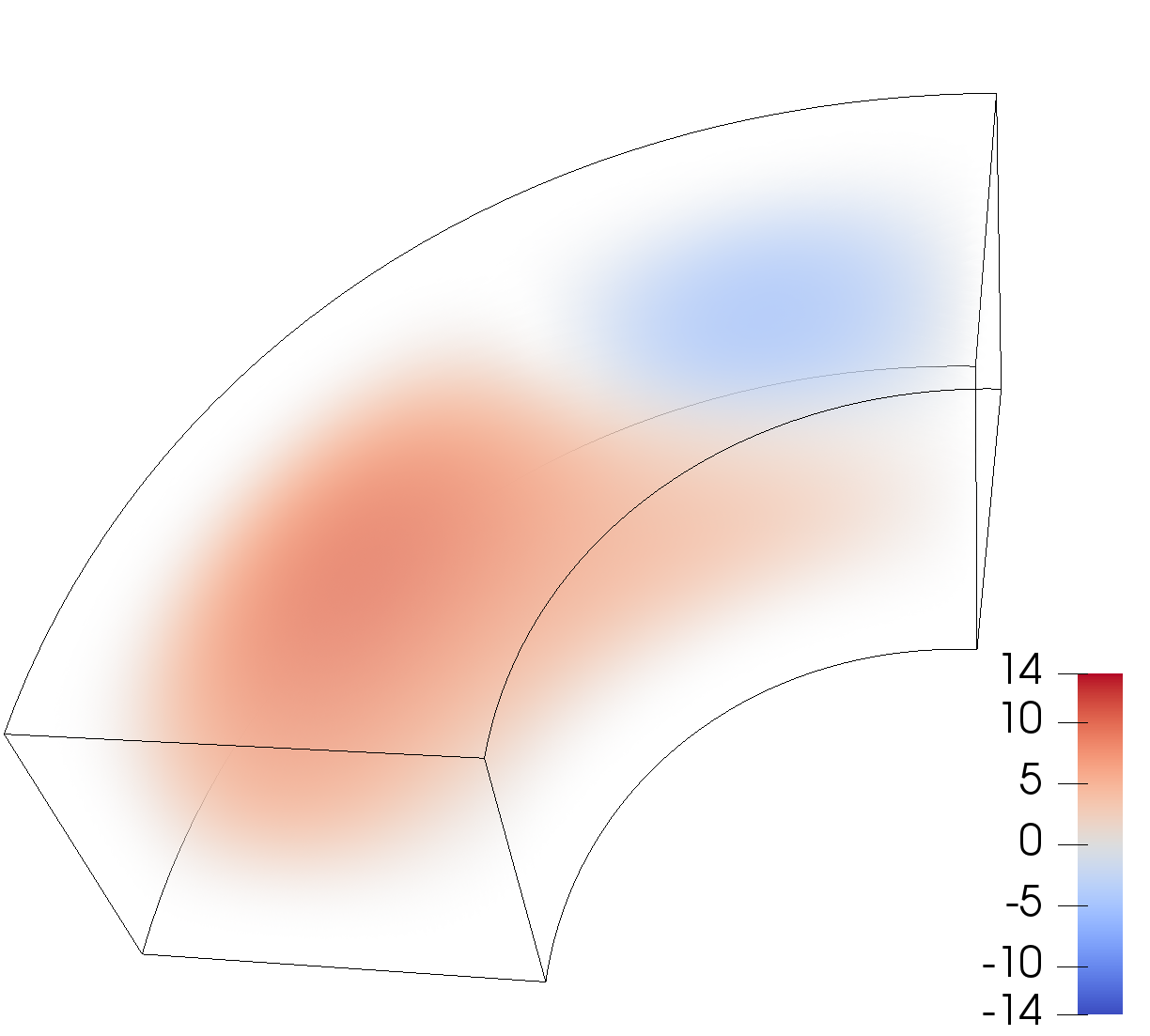}
   \caption{$\beta = 10^{-3}$} \label{figure:rank1_solution_1e-3}
  \end{subfigure}
  
  \begin{subfigure}[t]{0.3\textwidth}
  \centering
  \includegraphics[width=0.95\textwidth]{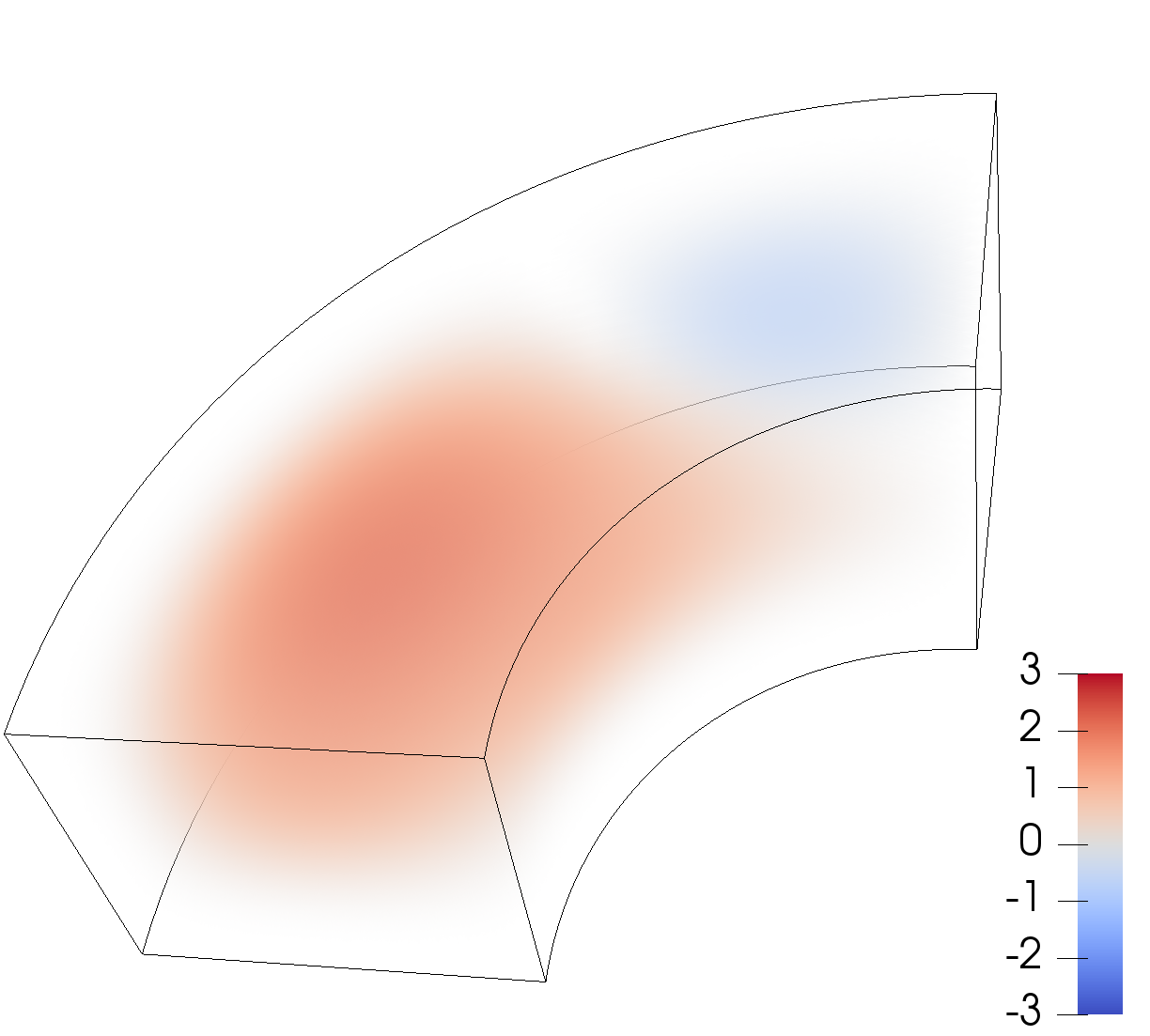}
  \caption{$\beta=10^{-2}$} \label{figure:rank1_solution_1e-2}
  \end{subfigure}%
    \begin{subfigure}[t]{0.3\textwidth}
   \centering
   \includegraphics[width=0.95\textwidth]{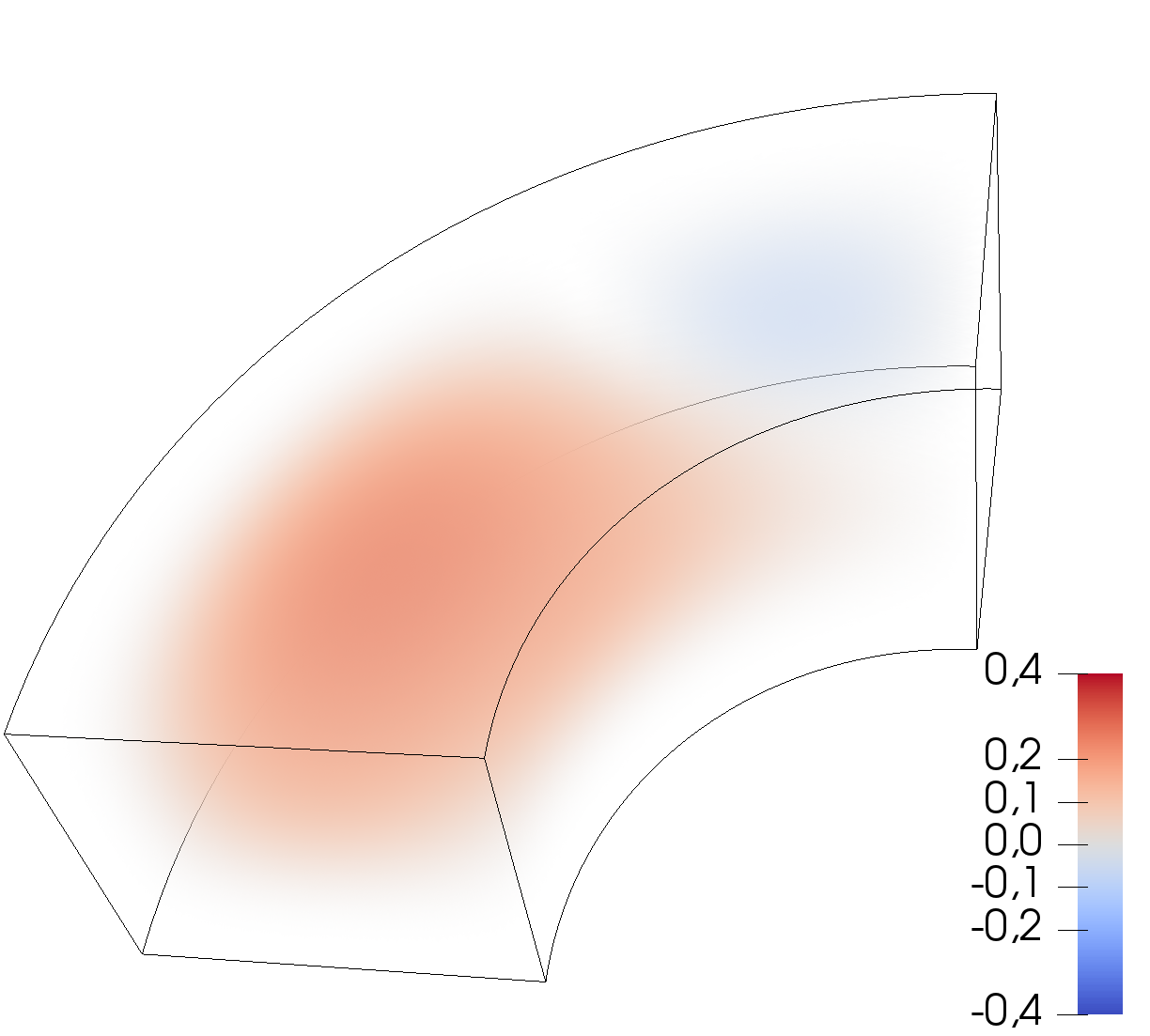}
   \caption{$\beta = 10^{-1}$} \label{figure:rank1_solution_1e-1}
  \end{subfigure}%
  \begin{subfigure}[t]{0.3\textwidth}
  \centering
  \includegraphics[width=0.95\textwidth]{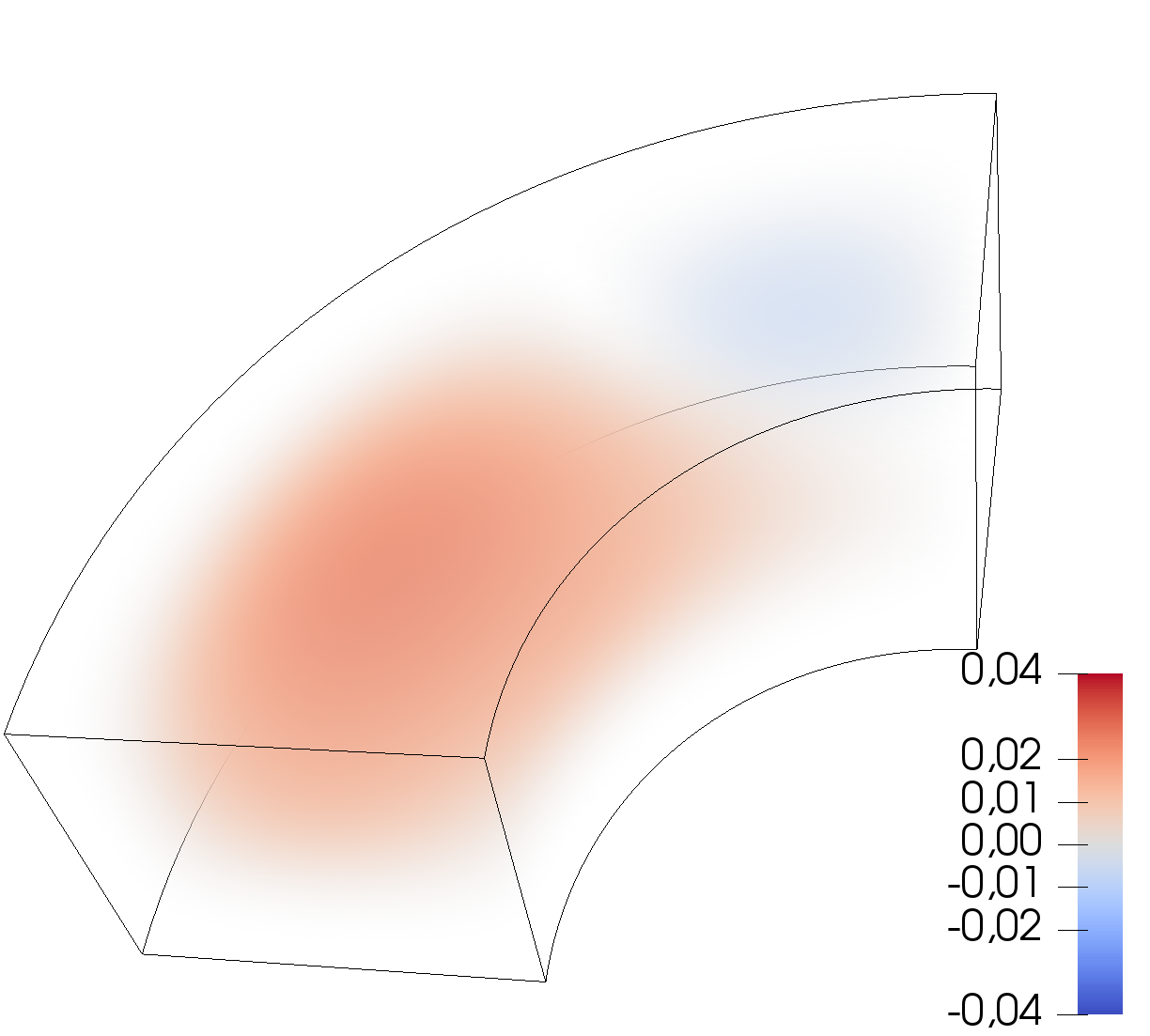}
  \caption{$\beta=1$} \label{figure:rank1_solution_1e0}
  \end{subfigure}

  \caption{Comparison of given state and solution for different control parameters at one time step on domain from Figure \ref{figure:geometry1}} \label{figure:rank1_solution}

\end{figure}

\begin{figure}
 \centering
 \includegraphics[width=0.9\textwidth]{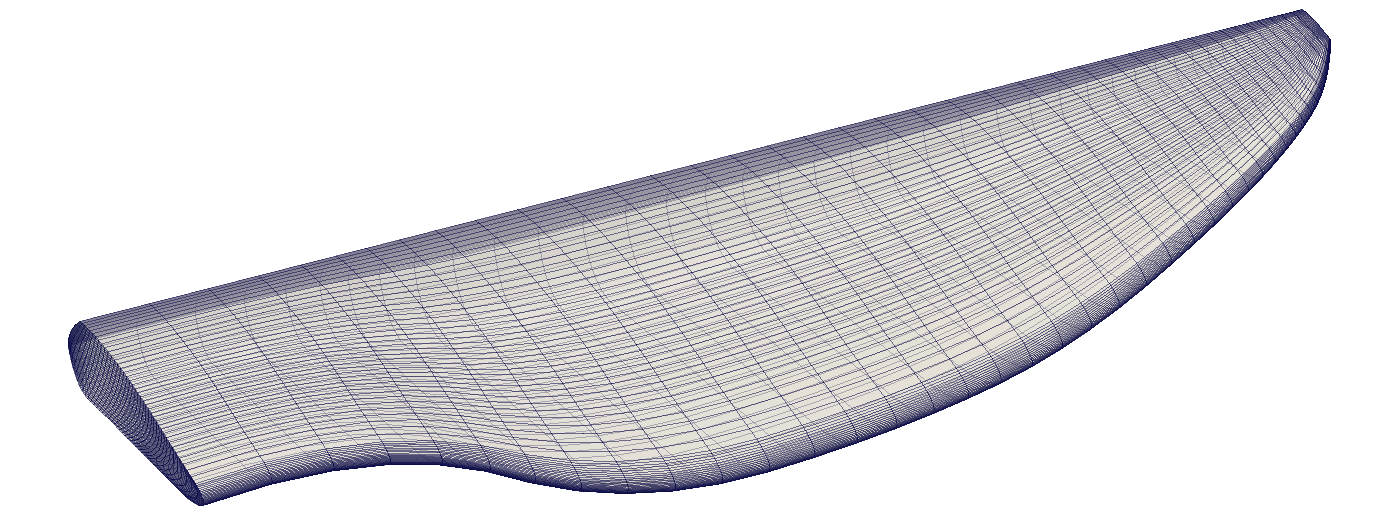}
 \caption{Freeform CAD model domain} \label{figure:geometryRotor}
\end{figure}

We now illustrate the performance of the Block AMEn method from Sec. \ref{section:LRoptimization} on two optimal control examples. We first regard the rank 1 domain from Figure \ref{figure:geometry1} before showing experimental results \tgreen{on the geometric model depicted in Figure \ref{figure:geometryRotor}}.

For the rank 1 geometry we used equidistant spatial knot insertion to refine the geometric representation and thus the solution space. \tgreen{We used discretizations with up to a maximum of 66 degrees of freedom per spatial direction and divided the time frame into 10 time steps. This discretization translates into a total of roughly 8 million degrees of freedom. Note that we will not pay any regard to the variation of the time discretization in this work. However, our experiments showed that increasing the number of time steps does not affect the number of iterative steps for most setups.}
The desired accuracy for \tblue{both the weight interpolation and the optimization was set to $10^{-5}$.}

In Figure \ref{figure:rank1opti2} we see the performance throughout different refinements for different control parameters $\beta$. Even for a high number of degrees of freedom the method converges after a small number of iterative steps as seen in Figure \ref{figure:rank1opti2iterations}. 

\tgreen{Using the Block AMEn method to solve the optimal control problem in a low rank format returns the solution in a low rank tensor train format too. Figure \ref{figure:rank1opti2ranks} illustrates the maximum TT-rank of the solution and even though they are quite large, the memory consumption of the solution is reduced drastically. Figure \ref{figure:rank1opti2memory} displays the storage requirements of the solution in relation to the full solution vector.}

Figure \ref{figure:rank1_solution} \tgreen{shows an exemplary result for the arbitrary desired state we used for our experiments.} The desired state in Figure \ref{figure:rank1_solution_yhat} was set as constant in time and Figures \ref{figure:rank1_solution_1e-4} - \ref{figure:rank1_solution_1e0} show snapshots of the same time step for different control parameters varying from $\beta=10^{-4}$ to $\beta=1$. As expected we see that the controlled state matches the desired state well for small control parameters and its magnitude decreases with higher control parameter.

The numerical values for the example from Figure \ref{figure:rank1_solution} are displayed in Figure \ref{figure:rank1opti1}. The method is robust with respect to the control parameter $\beta$. The objective function and the control behave as expected when the control parameter $\beta$ is changed and the method delivers a result within the desired accuracy after a small number of iterative steps.

\begin{figure}
 \centering
 \begin{subfigure}[t]{0.3\textwidth}
\centering
  \includegraphics[width=0.95\textwidth]{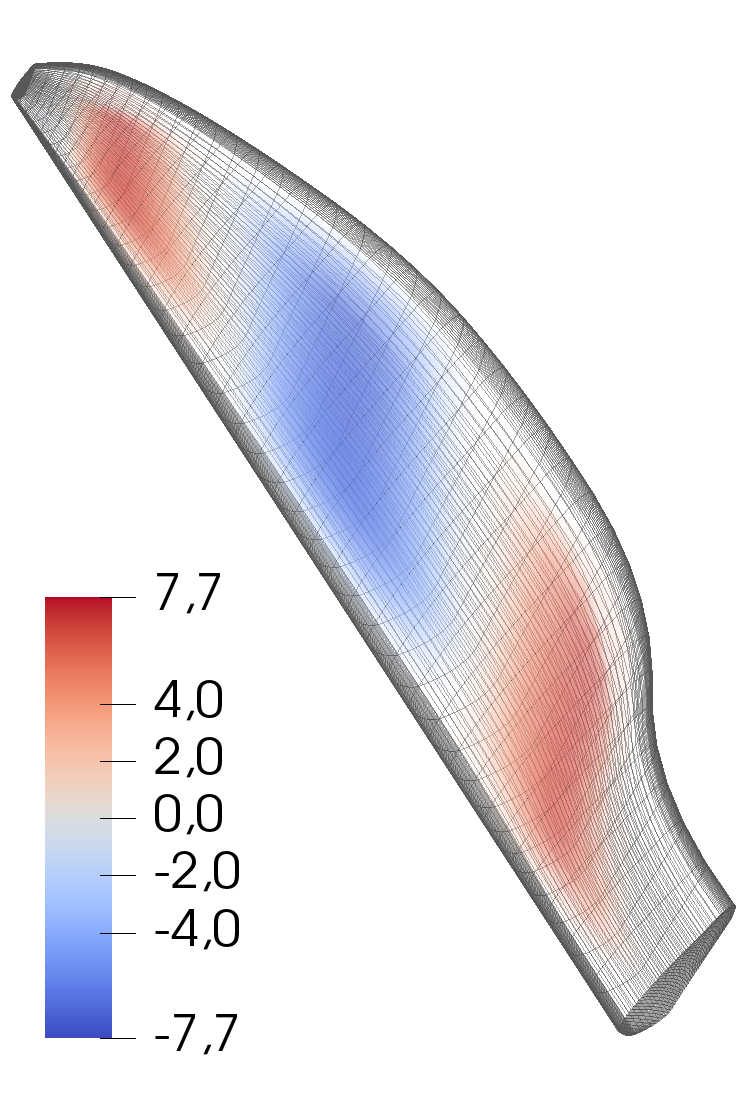} 
  \caption{Constant in time desired state $\hat{y}$} \label{figure:rotor_solution_yhat}
  \end{subfigure}%
  \begin{subfigure}[t]{0.3\textwidth}
   \centering
   \includegraphics[width=0.95\textwidth]{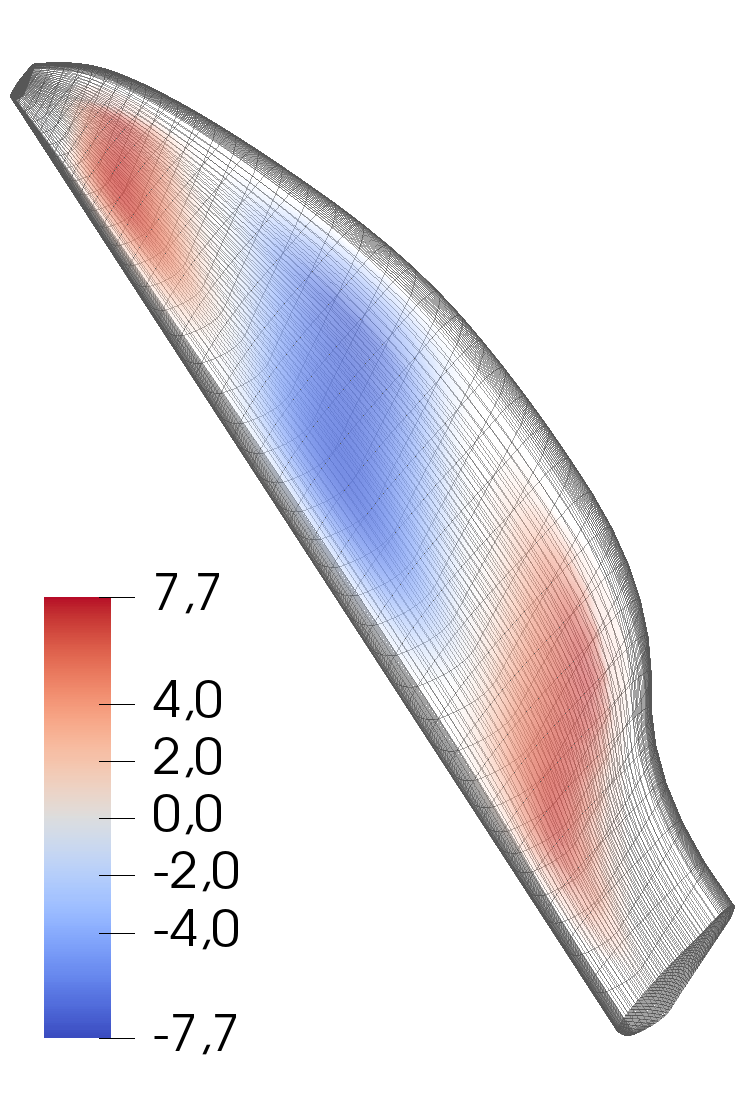}
   \caption{$\beta = 10^{-4}$} \label{figure:rotor_solution_1e-4}
   \end{subfigure}%
  \begin{subfigure}[t]{0.3\textwidth}
   \centering
   \includegraphics[width=0.95\textwidth]{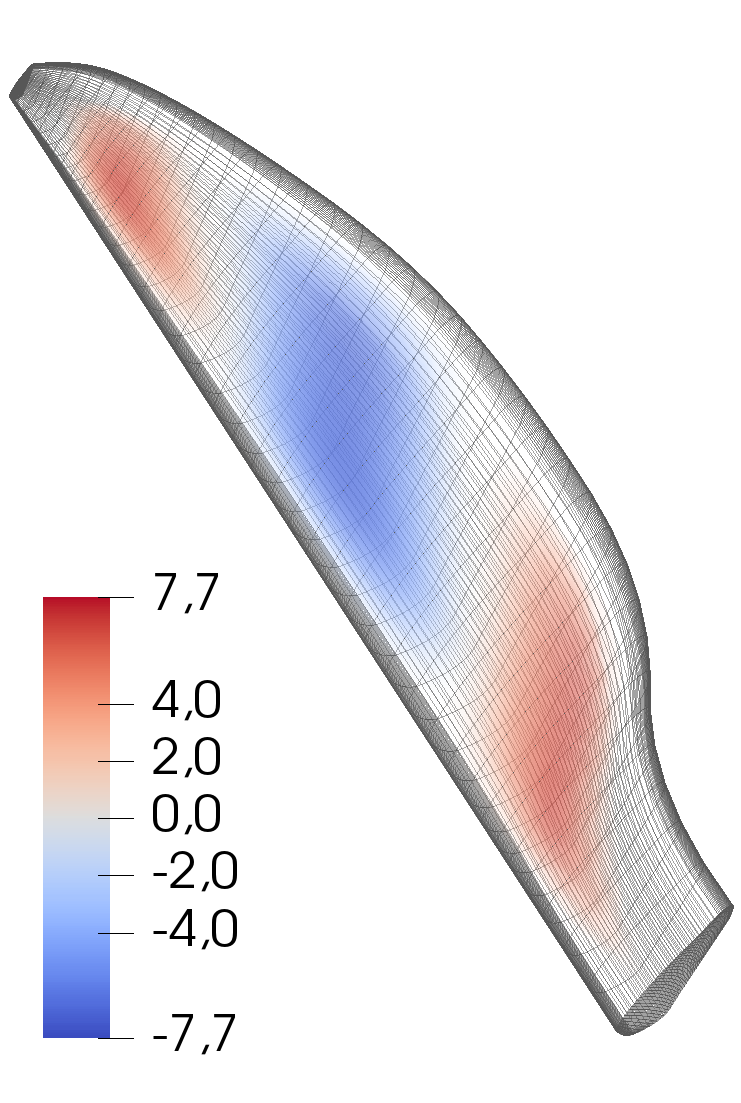}
   \caption{$\beta = 10^{-3}$}  \label{figure:rotor_solution_1e-3}
  \end{subfigure}
  
  \begin{subfigure}[t]{0.3\textwidth}
  \centering
  \includegraphics[width=0.95\textwidth]{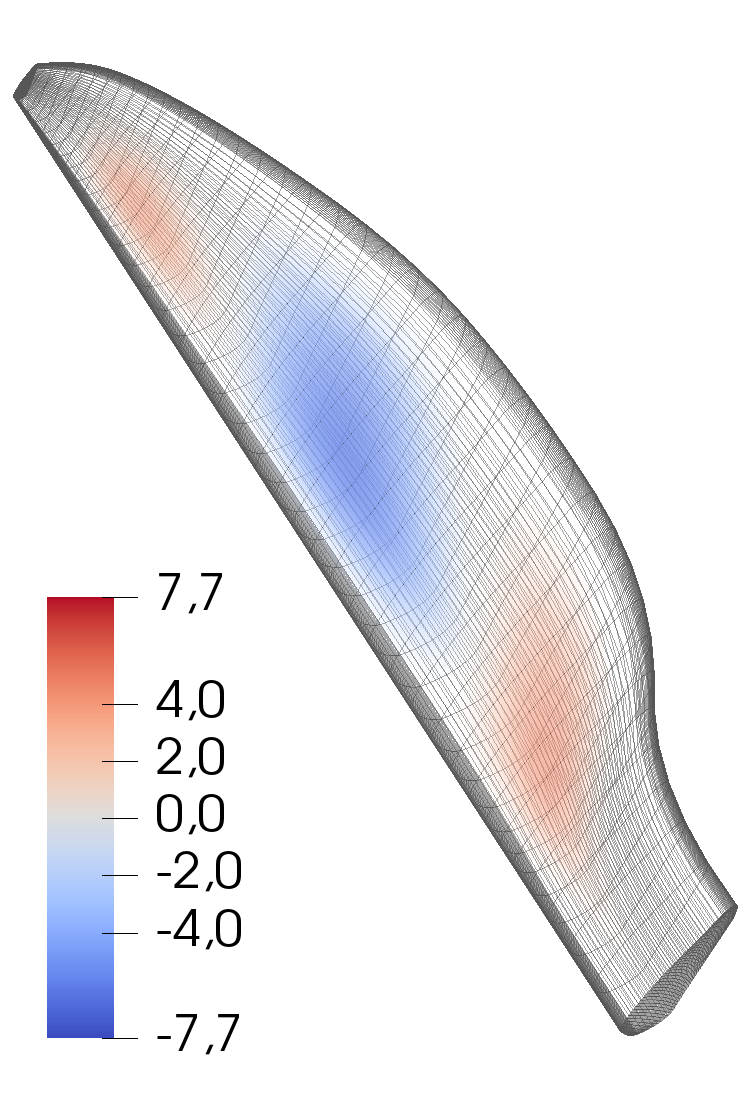} 
  \caption{$\beta=10^{-2}$} \label{figure:rotor_solution_1e-2}
  \end{subfigure}%
    \begin{subfigure}[t]{0.3\textwidth}
   \centering
   \includegraphics[width=0.95\textwidth]{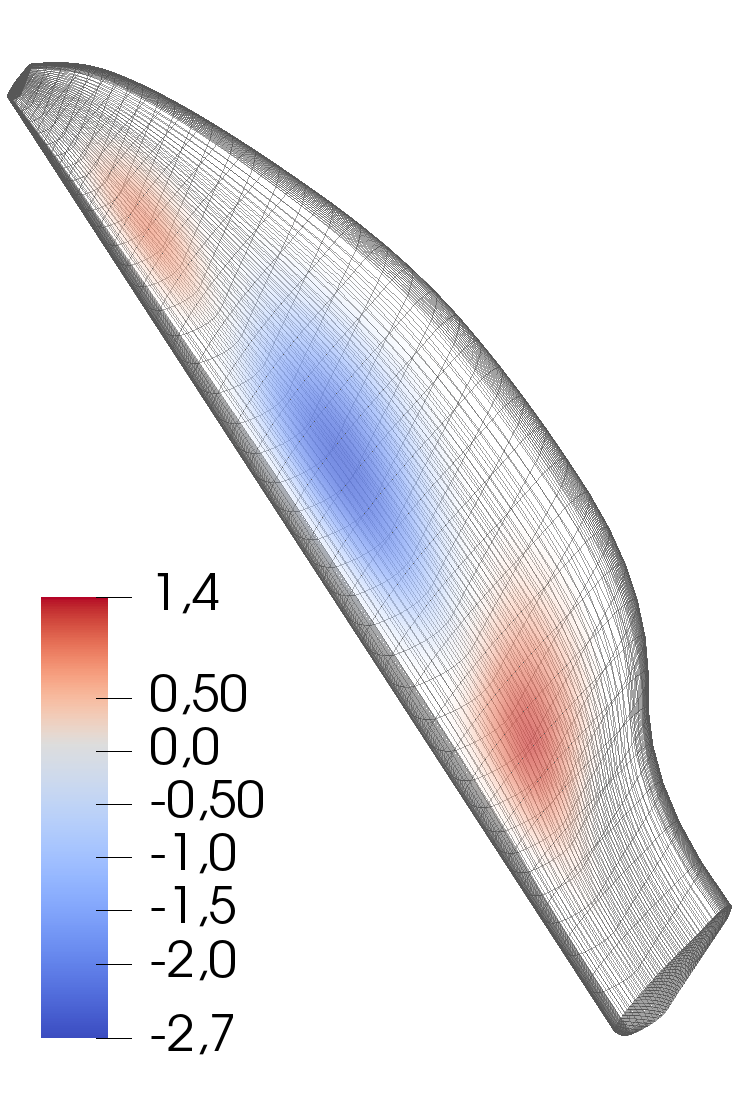}
   \caption{$\beta = 10^{-1}$} \label{figure:rotor_solution_1e-1}
  \end{subfigure}%
  \begin{subfigure}[t]{0.3\textwidth}
  \centering
  \includegraphics[width=0.95\textwidth]{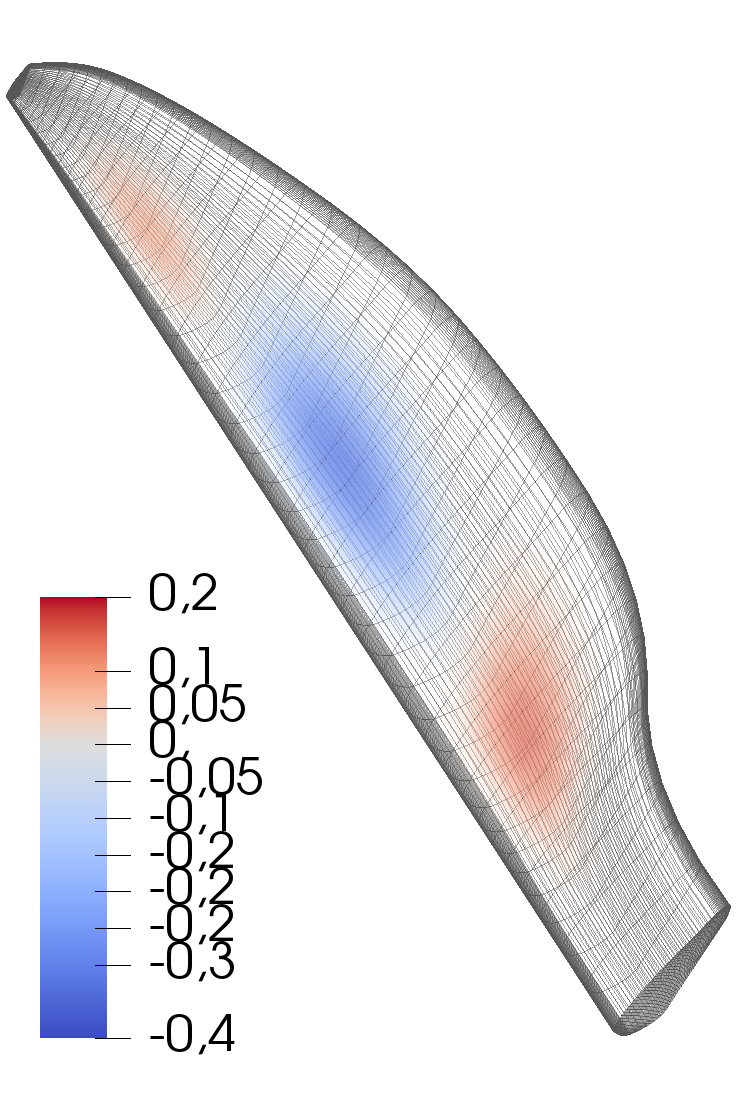}
  \caption{$\beta=1$} \label{figure:rotor_solution_1e0}
  \end{subfigure}

  \caption{Comparison of given state and solution for different control parameters at one time step on domain in Figure \ref{figure:geometryRotor}} \label{figure:rotor_solution}

\end{figure}

\tgreen{Additionally we tested the low rank optimization scheme on a large scale geometry inspired by the shape of a wind turbine rotor blade as depicted in Figure \ref{figure:geometryRotor}. The 3D solid NURBS model was designed by freeform surface modeling in the commercial CAD software Rhino 6.0. The low rank assembly step detected a rank profile as displayed in Table \ref{table:ranks3} for a truncation tolerance of $10^{-5}$.}

\begin{table}
 \centering
 \begin{tabular}{|c|c|c|c|c|c|c|c|c|}
 \hline
 tolerance & TT-ranks & $q_{1,1}$ & $q_{1,2}$ & $q_{1,3}$ & $q_{2,2}$ & $q_{2,3}$ & $q_{3,3}$ & $\omega$\\
 \hline
 $10^{-5}$ & $R_1$ & 11 & 12 & 5 & 11 & 3 & 4 & 4\\
 & $R_2$ & 2 & 3 &  1 & 2 & 1 & 2 & 1\\ 
\hline
 \end{tabular}
\caption{Ranks for the weight function tensor approximation of the CAD model in Figure \ref{figure:rotor_solution}} \label{table:ranks3}
\end{table}

\tgreen{Again, we set a fixed number of 10 time steps and an arbitrary desired state constant in time. The desired state used for the experiments is illustrated in Figure \ref{figure:rotor_solution_yhat}. Figures \ref{figure:rotor_solution_1e-4} - \ref{figure:rotor_solution_1e0} show a time snapshot of the experiment for different control parameters $\beta$. }

\tgreen{Even though this geometric model has a high rank profile, our scheme performes very well as seen in Figure \ref{figure:rotor_opti}.}
\tgreen{For larger control parameters $\beta$ the rank of the solution and the number of iterative steps are robust and stay almost constant for different levels of discretization as illustrated in Figures \ref{figure:rotor_opti_iterations} and \ref{figure:rotor_opti_ranks}. The number of iterations and the ranks increase only for very small control parameters. But even for the smallest control parameter with a high solution rank we reach a significant reduction in memory consumption comparing the low rank solution with the full solution vector as displayed in Figure \ref{figure:rotor_opti_memory}. }

\input{exp_optimization_rotor.tex}

\section{Conclusion}
In this paper, we combined the low rank method presented by Mantzaflaris et al. with Tensor Train calculations to obtain a powerful method for solving large equation systems arising from IGA-discretized PDEs \tgreen{and successfully applied the developed scheme to efficiently solve large PDE-constrained optimal control problems.} 

We can reduce the storage requirements and calculation time for the \tgreen{mass and} stiffness matrix assembly drastically by finding low rank approximations and splitting the matrices into a Kronecker product of smaller matrices. Our scheme finds low rank approximations for given desired accuracies without any prior knowlede about the geometry. \tgreen{We can exploit the} resulting low rank structures, keeping the memory consumption low throughout further computations. The iterative Block AMEn method allows us to solve large systems like a PDE-constrained optimal control problem without assembling the whole equation system. In combination with this iterative method the low rank format gives a great advantage and we can solve very large systems within a reasonably short time.

\tgreen{Various numerical experiments showed the high potential of the method. However, there might be even further efficiency gains if we} \tblue{find a suitable preconditioner for the reduced linear systems \eqref{eq:KKT_red} in the block AMEn method.}

\section*{Acknowledgment} The authors would like to thank Angelos Mantzaflaris for his helpful insights.
\bibliographystyle{siam}
\bibliography{literature}
\end{document}